\documentclass[preprint,10pt]{elsarticle}

\usepackage{amssymb}
\journal{Nonlinear Analysis-Real World Applications}

\usepackage{amssymb,amsmath,amsfonts,enumerate,bm}

\newtheorem{theorem}{Theorem}[section]
\newtheorem{lemma}[theorem]{Lemma}
\newtheorem{corollary}[theorem]{Corollary}
\newtheorem{proposition}[theorem]{Proposition}

\numberwithin{equation}{section}

\def\bproof{\noindent{\it Proof. }}
\def\eproof{\null\hfill {$\Box$}\bigskip}

\newcommand{\Er}{\mathbb{R}}

\newcommand{\en}{\mathbb{N}}

\newcommand{\Ers}{{\mathbb{R}}_{\text{sym}}^{3\times 3}}
\newcommand{\Erdev}{{\mathbb{R}}_{\text{dev}}^{3\times 3}}

\newcommand{\cK}{C^{\text{Korn}}}

\newcommand{\norm}[2][{}]{\lVert#2\rVert_{{#1}}}

\newcommand{\abs}[2][{}]{\lvert#2\rvert_{#1}}

\renewcommand{\tilde}{\widetilde}

\newcommand{\dive}[1]{\operatorname{div}#1}

\def\bfA{\bm{\mathrm A}}
\def\bfB{\bm{\mathrm{B}}}
\def\bfE{\bm{\mathrm{E}}}

\def\bfI{{\bm{\mathrm I}}}

\def\bfn{{\bm{\mathrm n}}}

\def\bfQ{{\bm{\mathrm Q}}}
\def\bfL{{\bm{\mathrm L}}}
\def\rmQ{\mathrm Q}
\def\bfsigma{{\bm{\mathrm{\sigma}}}}

\def\og{\leavevmode\raise.3ex\hbox{$\scriptscriptstyle\langle\!\langle$~}}
\def\fg{\leavevmode\raise.3ex\hbox{~$\!\scriptscriptstyle\,\rangle\!\rangle$}}

\def\C{{\mathrm{C}}}

\def\H{{\mathrm{H}}}
\def\L{{\mathrm{L}}}
\def\V{{\mathrm{V}}}
\def\X{{\mathrm{X}}}

\def\W{{\mathrm{W}}}
\def\ee{\bm\varepsilon}
\def\ze{z}
\def\zeo{z^{0}}
\def\dze{\dot z}

\def\tra{\mathsf{T}}
\def\eqldef{\overset{\text{\tiny def}}{=}}

\def\dd{\;\!\mathrm{d}}   %% d for integrals 

\def\calG{\mathcal G}

\def\calL{\mathcal L}
\def\calQ{\mathcal Q}

\def\calZ{\mathcal Z}

\begin{document}

\begin{frontmatter}



\title{Existence result for a class of generalized standard materials 
with thermomechanical coupling
{\footnote{This work was partially written when A.P. was employed by WIAS.
A.P. was supported  by the Deutsche Forschungsgemeinschaft 
through the projet C18 ``Analysis and numerics 
of multidimensional models for elastic phase transformation
in a shape-memory alloys'' of the Research Center {\sc Matheon}.
Moreover, L.P. gratefully acknowledges the hospitality of WIAS.}}}


\author[zob]{Laetitia Paoli\corref{cor1}}
\address[zob]{Universit\'e de Lyon, LaMUSE,
23 rue Paul Michelon\\ F-42023 Saint-Etienne Cedex 02, France\fnref{trs}}
\ead{laetitia.paoli@univ-st-etienne.fr}
\cortext[cor1]{Corresponding author: tel: +33(0)4.77.48.51.12,
fax: +33(0)4.77.48.51.53}

\author[boz]{Adrien Petrov}
\address[boz]{
Institut Camille Jordan UMR5208 and 
INSA de Lyon,
20 Avenue A. Einstein\\
F-69621 Villeurbanne Cedex,
France\fnref{trs}}
\ead{apetrov@math.univ-lyon1.fr}

%%%%%%%%%%%%%%%%%%%%%%%%%%%%%%%%%%%%%%%%%%%%%%%%%%%%%%%%%%%%%%%%%%%%%%%%%%%%%%%%

\begin{abstract}
This paper deals with the study of a three-dimensional model of 
thermomechanical coupling
for viscous solids exhibiting hysteresis effects. This model is written 
in accordance with the formalism of generalized standard materials
and it is composed of the momentum equilibrium equation
combined with the flow rule, which describes some stress-strain dependance, 
coupled to the heat-transfer equation. More precisely, 
the coupling terms  are linear with respect to the temperature and the displacement and non linear with respect to the internal variable.
% allowing to consider a large class of applications 
%including visco-elasto-plastic or shape-memory materials undergoing 
%thermal expansion. 
The main mathematical difficulty  lies in the fact that the natural framework for
the right-hand side of the heat equation is the space of  \(\L^1\) functions.
A local existence result for this thermodynamically 
consistent problem is obtained by using a fixed-point argument. 
Then the solutions are proved to be physically admissible  and global 
existence  is discussed under some additional assumptions
on the data.

%% Keywords are optional
\begin{keyword}
Existence result, generalized standard materials, differential inclusion, heat equation.
\MSC[2008] 35A01, 35K55, 74F05, 35Q80, 74C10, 74N30.
\end{keyword}
\end{abstract}

\end{frontmatter}

\section{Description of the problem}
\label{sec:description-problem}

Motivated by the study of visco-elasto-plastic materials and also
Shape-Memory Alloys (SMA), we consider in this paper a thermomechanical 
coupling for a class of Generalized Standard Materials (GSM) exhibiting 
hysteresis effects. More precisely, in the framework of GSM due to Halphen 
and Nguyen
(see \cite{HalNgu75SMSG}) the mechanical behavior of the material is 
described by the   {\emph{momentum equilibrium equation}} 
combined with a constitutive law ({\emph{flow rule}}) and the unknowns 
are the \emph{displacement} \(u\) and an \emph{internal variable}  \(\ze\) 
which allows to take into account some dissipation at the microscopic level. 
Indeed, plasticity and phase transitions are inelastic processes which involve 
some loss of energy, transformed into heat. Thus it is necessary to take 
into account the thermal process in the description of the problem.

The model considered here is based on the Helmholtz free energy 
\(W(\ee,\ze,\theta)\), depending on the 
\emph{infinitesimal strain tensor}  \(\ee=\ee(u)
\eqldef \tfrac12(\nabla u{+}\nabla
u^{\tra})\) for the \emph{displacement} \(u\),  
the  \emph{internal variable}  \(\ze\)
and the \emph{temperature} \( \theta\).
Here \((\cdot)^{\tra}\) denotes the transpose of a tensor.
We assume that $W$ can be decomposed as follows 
\begin{equation}
\label{eq:free_energy}
W(\ee,\ze,\theta)\eqldef
W^{\textrm{mech}}(\ee,\ze) - W^{\theta} (\theta)
+\theta W^{\textrm{coup}}(\ee,\ze),
\end{equation}
which ensures that entropy separates 
the thermal and mechanical variables (see \eqref{eq:s}). Let us emphasize  
that the last term in the right hand side of
\eqref{eq:free_energy} allows for coupling effects between 
the temperature and both the displacement and the internal variable.
We make the assumption of small deformations. The momentum equilibrium 
equation and the flow rule are given by
\begin{subequations}
\label{eq:ent_eqED}
\begin{align}
&-\dive(\bfsigma^{\textrm{\bf el}}
{+}\bfA\dot\ee)=f,
\label{eq:ent_eq_1ED}\\&
\label{eq:ent_eq_2ED}
\partial\Psi(\dze)+\bfB\dze+\bfsigma^{\textrm{\bf inel}}
\ni 0, 
\end{align}
\end{subequations}
where
$f$ is a given loading, 
\(\bfsigma^{\textrm{\bf el}}\eqldef\partial_{\ee} 
W(\ee,\ze, \theta)\),
\(\bfsigma^{\textrm{\bf inel}}
\eqldef\partial_{\ze} W(\ee,\ze, \theta) \),  \(\bfA\) and \(\bfB\) 
are two viscosity tensors 
and \(\Psi\) is the dissipation potential. As it is
common in modeling hysteresis effects in mechanics, 
we assume that \(\Psi\) is convex, positively homogeneous of
degree 1 and  $0 \in \partial \Psi (0)$ which ensures that 
\(\bfsigma^{\textrm{\bf inel}} {.} \dze \le 0 \).

Then the specific \emph{entropy} is defined by the Gibb's relation
\begin{equation} \label{eq:s}
s\eqldef -\partial_{\theta} W(\ee,z,\theta) 
=\partial_{\theta} W^{\theta} (\theta)-W^{\textrm{coup}} (\ee,z)  ,
\end{equation}
and the \emph{entropy equation} 
\begin{equation} \label{eq:entropy}
\theta\dot s - \dive ( \kappa \nabla \theta )=
\bfA \dot \ee{:} \dot \ee +
\bfB  \dot z{.}\dot z+\Psi(\dot z),
\end{equation}
gives some balance between the \emph{heat flux} 
$j = - \kappa \nabla \theta$, 
where \(\kappa\) is  the \emph{heat conductivity},  and the 
\emph{dissipation rate} $\xi\eqldef\bfA \dot \ee{:} \dot \ee +
\bfB  \dot z{.}\dot z+\Psi(\dot z) \ge 0$.
If the system is thermally isolated and $\theta>0$, we have 
\begin{equation*}
\int_{\Omega}\dot s\dd x=
\int_{\Omega}\tfrac{\dive(\kappa \nabla\theta)}{\theta}\dd x+
\int_{\Omega}\tfrac{\xi}{\theta}\dd x=
\int_{\Omega}\tfrac{\kappa \nabla\theta \cdot \nabla 
\theta}{\theta^2}\dd x+
\int_{\Omega}\tfrac{\xi}{\theta}\dd x \ge 0,
\end{equation*}
which guarantees that the second law of thermodynamics is satisfied. 
Furthermore, let 
\begin{equation*}
W^{\textrm{in}}(\ee,\ze ,\theta) \eqldef W (\ee,\ze ,\theta) + \theta s
\end{equation*}
be the \emph{internal energy}. By using the chain rule and  
\eqref{eq:ent_eqED}--\eqref{eq:entropy}, we obtain 
\begin{equation*}
\int_{\Omega} \dot W^{\textrm{in}}(\ee, \ze,\theta)\dd x
=
\int_{\Omega} f {\cdot}\dot u\dd x
 + \int_{\partial\Omega} \kappa\nabla\theta{\cdot}\bfn \dd x,
\end{equation*}
which gives the total energy balance in terms of the internal energy, 
the power of external load and heat.
Hence the model considered here is thermodynamically consistent.

We assume in the sequel that
\begin{subequations} \label{eq:energy}
\begin{align}
&W^{\textrm{mech}}(\ee,\ze)\eqldef
\tfrac12\bfE(\ee{-}\ee^{\textrm{\bf inel}}){:}(\ee{-}\ee^{\textrm{\bf inel}})+
\tfrac{\alpha}2\abs{\nabla\ze}^2+H_1(\ze),  \, \ee^{\textrm{\bf inel}} 
\eqldef\bfQ\ze, 
\label{eq:energy_1}\\&
W^{\theta}\eqldef c (\theta {\rm ln}(\theta){-}\theta), 
\label{eq:energy_2}\\&
\label{eq:energy_3}
W^{\textrm{coup}}(\ee,\ze)\eqldef\beta\bfI{:}\ee+H_2(\ze),
\end{align}
\end{subequations}
where  
\(c\) is the \emph{heat capacity}, 
\(\beta\geq 0\)
is the \emph{isotropic thermal expansion coefficient},
\(\bfI\) is the identity matrix, \(\alpha\geq 0\) is a coefficient 
that measures some non local
interaction effects for the internal variable \(\ze\),
\(\bfE\) is the \emph{elasticity tensor},
\(H_i\), \(i=1,2\), are two \emph{hardening functionals} and 
\(\bfQ\) is an affine mapping from a finite dimensional real vector 
space \(\calZ\) to \(\Ers\). More precisely, \(\bfQ\) is decomposed as follows
\begin{equation*}
\forall \ze\in \calZ:\ \bfQ \ze\eqldef \tilde\bfQ z+\rmQ,
\end{equation*}
with $\tilde \bfQ \in {\mathcal L}(\calZ, \Ers)$ and $\rmQ \in \Ers$.
We observe that by inserting 
on the one hand \eqref{eq:energy_1} and \eqref{eq:energy_3} into 
\eqref{eq:ent_eqED} and by carrying 
on the other hand \eqref{eq:energy_2}, 
\eqref{eq:energy_3} and \eqref{eq:s} into \eqref{eq:entropy}, we obtain
\begin{subequations}
\label{eq:ent_eq}
\begin{align}
&
-\dive(\bfE(\ee(u){-}\bfQ\ze){+}\beta\theta\bfI{+}\bfA\ee(\dot u))
=f,\label{eq:ent_eq_1}\\&
\partial\Psi(\dze)+\bfB\dze-\tilde\bfQ^{\tra}\bfE(\ee(u){-}\bfQ\ze)
+\partial_{\ze} H_1(\ze)+
\theta\partial_{\ze}H_2(\ze)-\alpha\Delta\ze\ni 0,\label{eq:ent_eq_2}\\&
\label{eq:ent_eq_3}
c\dot\theta{-}\dive(\kappa\nabla\theta)
=\bfA\ee(\dot u){:}\ee(\dot u)
+\theta(\beta\bfI{:}\ee(\dot u){+}\partial_{\ze}H_2(\ze){.}\dze)
+ \bfB\dze{.}\dze +\Psi(\dze),
\end{align}
\end{subequations}
together with boundary conditions
\begin{equation}
\label{eq:boun_cond}
u=0,\quad
\alpha \nabla \ze{\cdot}\bfn=0,\quad
\kappa\nabla\theta{\cdot}\bfn=0 \quad\text{ on }\quad
\partial\Omega\times [0,\tau),
\end{equation}
and initial conditions
\begin{equation}
\label{eq:init_cond}
u(\cdot,0)=u^0,\quad
\ze(\cdot,0)=\zeo,\quad
\theta(\cdot,0)=\theta^0\quad\text{ in }\quad
\Omega.
\end{equation}
Here $\Omega \subset \Er^3$ is a reference configuration and \(\bfn\) 
denotes the outward normal to the boundary \(\partial\Omega\)
of \(\Omega\). As usual, \((\,\dot{}\,)\), \(\partial_{\ze}^i\)
and \(\partial\) denote the time derivative
\(\tfrac{\partial}{\partial t}\), the \(i\)-th derivative  with respect 
to \(\ze\)
and the subdifferential in the
sense of convex analysis (see \cite{Brez73OMMS}),
respectively. 
Moreover \(\ee_1{:}\ee_2\) and \(\ze_1{.}\ze_2\) denote the
inner product of \(\ee_1\) and \(\ee_2\) in the  space of symmetric 
\(3{\times}3\) tensors \(\Ers\) and \(\ze_1\)
and \(\ze_2\) in the finite dimensional real vector space \(\calZ\). 

The increasing interest in smart materials for industrial applications 
has deeply stimulated the study of such models 
in engineering as well as in mathematical literature 
during the last decade. If 
the coupling with heat equation \eqref{eq:ent_eq_3} is ignored (for instance,  
if the characteristic dimension of the material is small in at least one 
direction, the temperature can be considered as a data), the problem 
\eqref{eq:ent_eq_1}--\eqref{eq:ent_eq_2} together with 
\eqref{eq:boun_cond}--\eqref{eq:init_cond} is nowadays quite well understood; 
existence results can be obtained either by using classical methods for 
maximal monotone operators (see \cite{AlbChe04QVMH}) 
or more specific techniques 
for rate-independent processes when the viscosity tensors vanish 
(see \cite{MieThe04RIHM,Miel05ERIS}). 
On the contrary, if the temperature is considered as 
an unknown, the coupling with the thermal process, which is not 
rate-independent, does not allow to use the previous techniques and 
the problem becomes much more difficult. Indeed, the natural functional 
framework for the right-hand side of \eqref{eq:ent_eq_3} seems at a first 
glance to be $\L^1(0,\tau; \L^1(\Omega))$ since we usually expect the 
displacements to be in $\W^{1,2}(0,\tau; \W^{1,2}(\Omega))$. This 
difficulty has been overcome in a serie of recent papers by using 
the so-called enthalpy transformation. More precisely, assuming that 
the heat conductivity is a continuous function of $\theta$ such that 
\begin{equation} 
\label{eq:enthalpy}
\exists \gamma >1, \ \exists c^c >0, \ \forall \theta \ge 0 : \ c(\theta) 
\ge c^c (1{+}\theta)^{\gamma -1},
\end{equation}
a new unknown, the \emph{enthalpy}, is defined by
\begin{equation*}
\vartheta \eqldef \int_0^{\theta} c(s) \dd s,
\end{equation*}
and the heat equation is replaced by the enthalpy equation:
\begin{equation*}
\dot\vartheta-\dive\bigl(\tfrac{\kappa}{c(\zeta(\vartheta))} 
\nabla\vartheta\bigr)
=\bfA \ee(\dot u){:}\ee(\dot u)
+ \zeta(\vartheta)(\beta\bfI{:}\ee(\dot u){+}\partial_{\ze}H_2(\ze){.}\dze)
+\bfB \dze{.}\dze +\Psi(\dze),
\end{equation*}
with $\abs{\theta{=}\zeta(\vartheta)} \leq \bigl(\tfrac{\gamma}{c^c} 
\max (\vartheta, 0) \bigr)^{\frac1{\gamma}}$. Roughly speaking, this change 
of unknown weakens the coupling effects (the greater is $\gamma$, 
the weaker is the coupling effects) and allows to build a solution 
either by using a time-discretization (\cite{BarRou08TPSS, Roub10TRIV, 
BarRou11TVEP}) or by using a fixed-point argument 
(\cite{PaoPet11GESM,PaoPet11GETV,PaoPet11TMPV}). Unfortunately, assumption 
\eqref{eq:enthalpy} on the heat conductivity is not always satisfied and 
we will consider in this paper  the more standard case where $c$ is a 
function of $x$. In such a case, the enthalpy is simply $\vartheta\eqldef c(x) 
\theta$ and does not provide any help in the mathematical analysis of the 
system  \eqref{eq:ent_eq}--\eqref{eq:init_cond}. In other words, we have to 
manage directly with the original coupling 
\eqref{eq:ent_eq_1}, \eqref{eq:ent_eq_2} and \eqref{eq:ent_eq_3}. 
For this problem, we will prove an existence result by using 
a fixed-point argument. Since the right-hand side of \eqref{eq:ent_eq_3} 
behaves as a quadratic term with respect to $\theta$, we can not expect a 
global existence result without some smallness assumptions on the coupling 
parameters $\beta$ and $\partial_{\ze} H_2$. 

The paper is organized as follows. In Section \ref{sec:statement-result}, 
we introduce the assumptions on the data, and we present the main result 
(local existence result). Then Section \ref{sec:proof-theorem} is 
devoted to its proof. In Section  \ref{sec:furth-prop-solut}, we 
establish some further properties of the solution, namely we prove 
that the temperature remains positive and thus is physically admissible 
and that $u$, $z$ and $\theta$ satisfy some global energy estimate. 
Furthermore we investigate sufficient conditions to 
get a global solution. Finally, in Section  \ref{sec:examples}, we 
present some examples which fit our modelization.

%%%%%%%%%%%%%%%%%%%%%%%%%%%%%%%%%%%%%%%%%%%%%%%%%%%%%%%%%%%%%%%%%%%%%%%%%%%%%%%%

\section{Statement of the result}
\label{sec:statement-result}

We consider a reference configuration $\Omega \subset \Er^3$, which is 
a bounded domain such that $\partial \Omega \in \C^{2+ \rho}$ with $\rho>0$.
Let us begin this section by introducing some assumptions on the data 
as well as obvious consequences following from these
assumptions used later on in this work.

\begin{enumerate}[({A}--1)]
\itemsep0.1em

\item 
The \emph{dissipation potential} $\Psi$ is positively homogeneous of 
degree $1$, satisfies the triangle inequality and remains bounded on 
the unit ball of ${\mathcal Z}$, i.e., we have
\begin{subequations}
\label{eq:Psi} 
\begin{align}
\label{hom}
&\forall \gamma \ge 0, \ \forall \ze\in \calZ: \ 
\Psi(\gamma \ze)=\gamma\Psi(\ze),
\\
\label{eq:4}
&\forall {\ze}_1,{\ze}_2\in \calZ:\ \Psi({\ze}_1{+}
{\ze}_2)\leq \Psi({\ze}_1)+\Psi({\ze}_2),
\\
\label{eq:Psi.bdd}
&\exists C^\Psi>0,  \ \forall \ze\in \calZ:\ 
0 \le \Psi(\ze)\leq C^\Psi \abs{\ze}.
\end{align}
\end{subequations} 
It is clear that \eqref{eq:Psi} 
implies that $\Psi$ is convex and continuous. With \eqref{eq:Psi.bdd}, 
we can also check immediately that $0 \in \partial \Psi(0)$.

\item The \emph{hardening functionals} \(H_i\), \(i=1, 2\),
 belong to  $\C^2 (\calZ; \Er)$ and satisfy the following inequalities
\begin{subequations}
\label{eq:H}
\begin{align}
\label{eq:H1}
&\exists c^{H_1}, \tilde c^{H_1}\geq 0 ,  \ \forall\ze\in\calZ :\ 
H_1(\ze)\geq c^{H_1}{\abs{\ze}}^2 - \tilde c^{H_1},
\\
\label{eq:H3}
&
\exists C^{H_i}_{\ze\ze}>0, \ \forall \ze\in\calZ:\ 
{\abs{\partial_{\ze}^2 H_i(\ze)}}\leq C^{H_i}_{\ze\ze}.
\end{align}
\end{subequations}

\item The \emph{elasticity tensor} \(\bfE:\Omega\rightarrow \calL(\Ers;\Ers)\)
is a symmetric  positive definite operator such that
\begin{subequations}
\label{eq:47}
\begin{align}
\label{eq:3}
&
\exists c^{\bfE}>0, \ \forall \ee\in\L^2(\Omega; \Ers):\ 
c^{\bfE}\norm[\L^2(\Omega)]{\ee}^2
\leq
\int_{\Omega}\bfE \ee{:}\ee\dd x,\\&
\label{eq:46}
\forall i,j,k=1,2,3:\
\bfE(\cdot),
\tfrac{\partial\mathrm{E}_{i,j}(\cdot)}{\partial x_k}\in
\L^{\infty}(\Omega).
\end{align}
\end{subequations}

\item The \emph{viscosity tensors} \(\bfA\) and 
\(\bfB\) are symmetric  positive 
definite such that
\begin{subequations}
\label{eq:LM} 
\begin{align}
\label{eq:6}
&
\exists c^{\bfA}, C^{\bfA}>0, \ \forall \ee\in\Ers:\ 
c^{\bfA}\abs{\ee}^2
\leq
\bfA \ee{:} \ee\leq
C^{\bfA}\abs{\ee}^2,
\\
\label{eq:8}
&
\exists c^{\bfB}, C^{\bfB}>0, \ \forall \ze\in\calZ:\ 
c^{\bfB}\abs{\ze}^2\leq
\bfB \ze{.}\ze\leq
C^{\bfB}\abs{\ze}^2.
\end{align}
\end{subequations}

\item The \emph{inelastic strain} is given by 
$\ee^{\textrm{\bf inel}}\eqldef\bfQ\ze = \tilde \bfQ \ze + \rmQ$ with
\begin{equation} \label{eq:Q}
\tilde \bfQ \in  {\mathcal L}(\calZ, \Ers) \quad \textrm{ and }\quad
\rmQ \in \Ers.
\end{equation}

\item The \emph{external loading} \(f\) satisfies
\begin{equation} \label{eq:f}
f \in \H^1(0,T;\L^2(\Omega))\quad{\text{with}}\quad T>0.
\end{equation}

\item The \emph{heat capacity}
\(c : \Omega \rightarrow \Er\) and 
the \emph{conductivity} \(\kappa^c: \Omega \rightarrow \Ers\) 
satisfy the following inequalities
\begin{subequations}
\label{eq:36}
\begin{align}
\label{eq:40}
&\exists C^c, c^c>0:\
c^c\leq c(x)\leq C^c \ \text{ a.e. } \ x\in\Omega,\\&
\label{eq:41}
\exists c^{\kappa}>0, \ \forall v\in\Er^3:\ 
\kappa(x) v{.}v\geq 
c^{\kappa}\abs{v}^2 \ \text{ a.e. } \ x\in\Omega, \\&
\label{eq:41b}
\exists  C^{\kappa}>0:\ \abs{\kappa(x)} \leq 
C^{\kappa}\ \text{ a.e. }\ x\in\Omega.
\end{align}
\end{subequations}
\end{enumerate}
Finally, we assume that $\alpha \ge 0$ and either 
$\alpha >0$ and $c^{H_1}>0$ or $\alpha =0$ and $\partial_{\ze} H_2 \equiv 0$. 
Note that the boundary condition $\alpha \nabla\ze {\cdot}\bfn=0$
on $ \partial\Omega\times[0,\tau)$ will disappear if $\alpha =0$.
Let us observe also that
\eqref{eq:H3} leads to
\begin{equation}
\begin{aligned}
\label{eq:H2}
&\exists C^{H_i}_{\ze}>0, \ \forall \ze\in\calZ : \\&
{\abs{\partial_{\ze} H_i(\ze)}}\leq C^{H_i}_{\ze}(1{+}{\abs{\ze}}),\quad
{\abs{H_i(\ze)}} \le  C^{H_i}_{\ze}(1{+}{\abs{\ze}}^2).
\end{aligned}
\end{equation} 
We use later the following notations; 
\(\W_{\textrm{Dir}}^{m,r}(\Omega)\eqldef \{
\xi\in\W^{m,r}(\Omega):\ \xi=0\; \text{ on }\;  \partial
\Omega\}\) and \(\W_{\textrm{Neu}}^{m,r}(\Omega)\eqldef \{
\xi\in\W^{m,r}(\Omega):\ \nabla\xi{\cdot}\bfn=0\; \text{ on }\;  \partial
\Omega\}\) where \(m\geq 1\) and \(r\geq 2\) are two integers.

As usual  Korn's inequality will play a role in the
mathematical analysis developed in the next sections. 
We have assumed that $\partial \Omega$ is of class $\C^{2+\rho}$, 
so we have 
\begin{equation}
\label{eq:58} 
\exists \cK>0,\ \forall u\in \W^{1,2}_{\textrm{Dir}}(\Omega):\
\norm[\L^2(\Omega)]{\ee(u)}^2
\geq
\cK
\norm[\W^{1,2} (\Omega)]{u}^2,
\end{equation}
for further details on Korn's
inequality, the reader is referred to \cite{DuvLio76IMP,KonOle88BVPS}. 

As a starting point in the study of the problem 
\eqref{eq:ent_eq}--\eqref{eq:init_cond}, let us consider the heat equation
\begin{equation}
\label{eq:19}
c\dot\theta-\dive(\kappa\nabla\theta)=f^{\tilde \theta},
\end{equation}
with initial and boundary conditions
\begin{equation}
\label{eq:21}
\theta(\cdot,0)=\theta^0 \quad\text{in}\quad \Omega, \quad 
\kappa\nabla\theta{\cdot}\bfn=0
\quad\text{on}\quad\partial\Omega\times[0,\tau).
\end{equation}
If $f^{\tilde \theta} \in \L^2(0,\tau; \L^2(\Omega))$ and 
$\theta^0 \in \W^{1,2}_{\textrm{$\kappa$, Neu}} 
(\Omega)\eqldef \{ \xi \in \W^{1,2}( \Omega):\  
\kappa\nabla\xi{\cdot}\bfn=0 \; \text{ on }\;  \partial
\Omega \}$, this problem admits a unique solution $\theta \in 
\L^{\infty} (0,\tau;\W^{1,2}_{\textrm{Neu}} (\Omega))$
(see \cite{Evan10PDES}).

Now we recall existence and uniqueness results 
for the system composed by the momentum equilibrium 
equation and the flow rule \eqref{eq:ent_eq_1}--\eqref{eq:ent_eq_2}
when the temperature is a given data.
More precisely, let \(\tilde \theta\) be given 
in \(\L^q(0,\tau ;\L^p(\Omega))\).
We consider  the following problem: Find $u : 
[0,\tau] \to \Er^3$ and $\ze: [0, \tau] \to {\mathcal Z}$ such that 
\begin{subequations}
\label{eq:ent_eq_1_1}
\begin{align}
&-\dive(\bfE(\ee(u){-}\bfQ\ze){+}\beta \tilde \theta\bfI{+}\bfA\ee(\dot u))
=f,\label{eq:ent_eq_1_2}\\&
\partial\Psi(\dze)+\bfB\dze-\tilde\bfQ^{\tra}\bfE(\ee(u){-}\bfQ\ze)
+\partial_{\ze} H_1(\ze)+
\tilde \theta\partial_{\ze}H_2(\ze)-\alpha\Delta\ze\ni 0,\label{eq:ent_eq_2_2}
\end{align}
\end{subequations}
with boundary conditions
\begin{equation}
\label{eq:boun_cond2}
u=0,\quad
\alpha \nabla \ze{\cdot}\bfn=0 \quad\text{ on }\quad
\partial\Omega\times [0,\tau),
\end{equation}
and initial conditions
\begin{equation}
\label{eq:init_cond2}
u(\cdot,0)=u^0,\quad
\ze(\cdot,0)=\zeo
\quad\text{ in }\quad\Omega.
\end{equation}

\begin{proposition}
\label{sec:existence}
Let $\tau \in (0, T]$ and $\tilde \theta$ be given in 
$\L^q(0, \tau; \L^p(\Omega))$ with $p\in [4,6] $. 
Assume that \eqref{eq:Psi}, \eqref{eq:H}, \eqref{eq:47}, 
\eqref{eq:LM}, \eqref{eq:Q}, \eqref{eq:f}, 
$u^0 \in \W^{1,p}_{\textrm{Dir}}(\Omega)$ and $\zeo
\in \W^{2,p}_{\textrm{Neu}}(\Omega)$
if \(\alpha>0\) or $\zeo\in \L^p(\Omega)$ if \(\alpha=0\)
hold. Then the problem 
\eqref{eq:ent_eq_1_1}--\eqref{eq:init_cond2}
admits a unique solution 
\(u\in\W^{1,q}(0,\tau;\W^{1,p}_{\rm{Dir}}(\Omega))\) and 
\(\ze\in\L^{q/2} (0,\tau;\W^{2,p}_{\rm{Neu}}(\Omega)) 
\cap \C^0([0, \tau]; \W^{1,2}_{\rm{Neu}}(\Omega))
\cap\W^{1,q/2}(0,\tau;\L^{p}(\Omega))
\cap\W^{1,q}(0,\tau;\L^{p/2}(\Omega))\) if \(\alpha>0\) and
\(\ze\in\W^{1,q}(0,\tau;\L^p(\Omega))\) if \(\alpha=0\)
for any \(q>8\). Furthermore $\tilde \theta \mapsto (u,\ze)$ 
maps any bounded subset of $\L^q(0, \tau; \L^p(\Omega))$ into 
a bounded subset of $ \W^{1,q}(0,\tau;\W^{1,p}_{\rm{Dir}}(\Omega))
\times(\L^{q/2} (0,\tau;\W^{2,p}_{\rm{Neu}}(\Omega)) 
\ \cap\ \C^0([0, \tau]; \W^{1,2}_{\rm{Neu}}(\Omega))
\ \cap\ \W^{1,q/2}(0,\tau;\L^{p}(\Omega))
\ \cap\ \W^{1,q}(0,\tau;\\ \L^{p/2}(\Omega)))$ if \(\alpha>0\) or
into a bounded subset of
$\W^{1,q}(0,\tau;\W^{1,p}_{\rm{Dir}}(\Omega)){\times} 
\W^{1,q}(0,\tau;\\ \L^p(\Omega))$ if \(\alpha=0\). 
\end{proposition}

The key tool to prove existence, uniqueness and boundedness results for
\eqref{eq:ent_eq_1_1}--\eqref{eq:init_cond2} consists in 
interpreting this system of partial differential equations
 as an ordinary differential equation  in an appropriate Banach space. 
For the detailed proof, the reader is referred
to \cite[Thm.~4.1, Prop.~4.2, Lem.~4.4--4.5]{PaoPet11GESM} and 
\cite[Thm.~3.1, Prop.~3.2, Lem.~3.4--3.5]{PaoPet11TMPV} 
when $\alpha >0$ and to
\cite[Thm.~4.1]{PaoPet11GETV} when $\alpha =0$. 

So we may prove the existence of a solution  for the coupled problem 
\eqref{eq:ent_eq}--\eqref{eq:init_cond} by combining via a 
fixed-point argument the results of Proposition \ref{sec:existence} 
with the existence results for the heat equation with $f^{\tilde \theta}$ 
given by
\begin{equation*}
f^{\tilde \theta}\eqldef \bfA\ee(\dot u){:}\ee(\dot u)
+\tilde\theta(\beta\bfI{:}\ee(\dot u){+}\partial_{\ze}H_2(\ze){.}\dot\ze) 
+\bfB\dot\ze{.}\dot \ze +\Psi(\dot\ze).
\end{equation*}
We will obtain

\begin{theorem}
\label{thm:local_existence}
Assume that \eqref{eq:Psi}, \eqref{eq:H}, \eqref{eq:47}, 
\eqref{eq:LM},  \eqref{eq:Q}, \eqref{eq:f}, \eqref{eq:36},
\(\theta^0\in\W^{1,2}_{\kappa, \rm{Neu}}(\Omega)\), 
\(u^0\in\W^{1,4}_{\rm{Dir}}(\Omega)\)
and \(\zeo\in\W^{2,4}_{\rm{Neu}}(\Omega)\) if \(\alpha>0\) and
\(\zeo\in\L^4(\Omega)\) if \(\alpha=0\) hold. Then there exists 
\(\tau\in(0,T]\) such that the problem \eqref{eq:ent_eq}--\eqref{eq:init_cond}
admits a solution on \([0,\tau]\) such that
\(\theta\in\L^{\infty}(0,\tau;\W^{1,2}_{\kappa, \rm{Neu}}(\Omega))
\cap\C^0(0,\tau;\L^4(\Omega))\), 
 \(\dot\theta\in\L^2(0,\tau;\L^2(\Omega))\), 
\(u\in\W^{1,q}(0,\tau;\W^{1,4}_{\rm{Dir}}(\Omega))\) and  
\(\ze\in\L^{q/2} (0,\tau;\W^{2,4}_{\rm{Neu}}(\Omega))
\cap \C^0([0, \tau]; \W^{1,2}_{\rm{Neu}}(\Omega))
\cap\W^{1,q/2}(0,\tau;\L^{4}(\Omega))
\cap\W^{1,q}(0,\tau;\L^{2}(\Omega))\) when \(\alpha>0\),
\(\ze\in\W^{1,q}(0,\tau;\L^4(\Omega))\)
when \(\alpha=0\), for any $q>8$.
\end{theorem}

Next, recalling that the problem is thermodynamically consistent if 
$\theta>0$, we establish at Proposition \ref{temp_adm} that the 
solution obtained  in the previous theorem is physically admissible, 
i.e. remains positive whenever $\theta^0 \ge \bar \theta$ almost 
everywhere in $\Omega$, with $\bar \theta>0$. Finally, a 
global energy estimate is obtained in Proposition \ref{thm:global-estim} and 
sufficient conditions on $\beta$ and $\partial_{\ze} H_2$ are proposed
to get a global solution $(u,\ze,\theta)$ defined on $[0,T]$.

%%%%%%%%%%%%%%%%%%%%%%%%%%%%%%%%%%%%%%%%%%%%%%%%%%%%%%%%%%%%%%%%%%%%%%%%%%%%%%%%

\section{Proof of Theorem \ref{thm:local_existence}}
\label{sec:proof-theorem}

This section is dedicated to the proof of 
Theorem \ref{thm:local_existence}   
by using a fixed-point argument.
More precisely, for any given $\tilde \theta\in\C^0([0,\tau];\L^4(\Omega))$
with $\tau \in (0,T]$,
we consider \(f^{\tilde \theta}\eqldef \bfA\ee(\dot u){:}\ee(\dot u)
+\tilde\theta(\beta\bfI{:}\ee(\dot u){+}\partial_{\ze}H_2(\ze){.}\dot\ze) 
+\bfB\dot\ze{.}\dot \ze +\Psi(\dot\ze)\) where \((u,\ze)\)
is  the unique  solution of \eqref{eq:ent_eq_1_1}--\eqref{eq:init_cond2}.
Using Proposition \ref{sec:existence}, we obtain
$f^{\tilde \theta} 
\in \L^2(0,\tau; \L^2(\Omega))$ and thus
the heat-transfer equation 
\eqref{eq:19}--\eqref{eq:21} possesses a unique solution 
$\theta \in \L^{\infty}( 0, \tau; \W^{1,2}_{\textrm{Neu}}(\Omega))$ such that 
$\dot \theta \in \L^2(0,\tau;\L^2(\Omega))$.
This allows us to define a mapping
\begin{equation*}
\begin{aligned}
\Phi_{\tau}^{\tilde\theta,\theta}:\  
\C^0([0,\tau];\L^4(\Omega))&\rightarrow \C^0([0,\tau];\L^4(\Omega))\\
\tilde \theta &\mapsto \theta.
\end{aligned} 
\end{equation*}
Our aim consists in proving that this mapping 
satisfies the assumptions of Schauder's fixed point theorem for some positive 
$\tau \in (0, T]$.

Let us define the set \(\calQ_{\tau}\eqldef \Omega\times (0,\tau)\) with
\(\tau\in(0,T]\). In the sequel, the notations for the constants introduced 
in the
proofs are valid only in the proof.

\begin{proposition}
\label{prop:Phi1}
Let \(\tau\in (0,T]\). Assume that \eqref{eq:Psi}, \eqref{eq:H}, 
\eqref{eq:47}, \eqref{eq:LM}, \eqref{eq:Q}, \eqref{eq:f}, \eqref{eq:36}, 
\(\theta^0\in\W^{1,2}_{\kappa, \rm{Neu}}(\Omega)\), 
\(u^0\in\W^{1,4}_{\rm{Dir}}(\Omega)\)
and \(\zeo\in\W^{2,4}_{\rm{Neu}}(\Omega)\) if \(\alpha>0\) and
\(\zeo\in\L^4(\Omega)\) if \(\alpha=0\) hold.
Then
\(\Phi_{\tau}^{\tilde\theta,\theta}\) maps any bounded subset of 
\(\C^0([0,\tau];\L^4(\Omega))\) into a bounded relatively compact subset
of \(\C^0([0,\tau];\L^4(\Omega))\).
\end{proposition}

\bproof
We recall first existence, uniqueness and regularity results
for the heat-transfer equation.
More precisely, let consider the system \eqref{eq:19}--\eqref{eq:21}.
We assume that \eqref{eq:36} holds and that the initial temperature 
\(\theta^0\in\W^{1,2}_{\textrm{$\kappa$, Neu}}(\Omega)\) and 
$f^{\tilde \theta}\in\L^2(0,\tau;\L^2(\Omega))$.  By using Galerkin's 
method (see for instance \cite{Evan10PDES}), we may 
prove that this problem admits a unique
solution \(\theta\in
\L^{\infty} (0,\tau;\W^{1,2}_{\textrm{Neu}}(\Omega))\) with $\dot\theta \in
\L^2(0,\tau;\L^2(\Omega))$. 

Moreover we have the following a priori estimates
\begin{equation}
\begin{aligned}
\label{eq:2}
&\norm[\L^2(\Omega)]{\theta (\cdot,t )}^2 
+ \tfrac{2c^{\kappa}}{c^c}
\norm[\L^2(0,\tau;\L^2(\Omega))]{\nabla\theta}^2\\&
\leq
\tfrac1{c^c}
\bigl(C^c\norm[\L^2(\Omega)]{\theta^0}^2
{+}
\norm[\L^2(0,\tau; \L^2(\Omega))]{f^{\tilde \theta}}^2\bigr)
\exp\bigl(\tfrac{\tau}{c^c}\bigr)
\end{aligned}
\end{equation}
and 
\begin{equation}
\begin{aligned}
\label{eq:5}
&c^c\norm[\L^2(0,\tau;\L^2(\Omega))]{\dot \theta}^2+
c^{\kappa}\norm[\L^2(\Omega)]{\nabla\theta(\cdot,t)}^2\\&
\leq
C^{\kappa}\norm[\L^2(\Omega)]{\nabla\theta^0}^2+
\tfrac1{c^c}\norm[\L^2(0,\tau;\L^2(\Omega))]{f^{\tilde \theta}}^2
\end{aligned}
\end{equation}
for almost every $t \in [0, \tau]$.
Therefore adding \eqref{eq:2} and \eqref{eq:5}, we have
\begin{equation}
\begin{aligned}
\label{eq:7}
&c^c\norm[\L^2(0,\tau;\L^2(\Omega))]{\dot\theta}^2
+\min(1,c^{\kappa})\norm[\W^{1,2}(\Omega)]{\theta(\cdot,t)}^2+
\tfrac{2c^{\kappa}}{c^c}\norm[\L^2(0,\tau;\L^2(\Omega))]{\nabla\theta}^2
\\&\leq
\max\bigl(\tfrac{C^c}{c^c}\exp\bigl(\tfrac{\tau}{c^c}\bigr),C^{\kappa}\bigr)
\norm[\W^{1,2}(\Omega)]{\theta^0}^2+\tfrac1{c^c}
\bigl(\exp\bigl(\tfrac{\tau}{c^c}\bigr){+}1\bigr)
\norm[\L^2(0,\tau;\L^2(\Omega))]{f^{\tilde \theta}}^2
\end{aligned}
\end{equation}
for almost every $t \in [0, \tau]$.
We introduce now the following functional space
\begin{equation*}
\V((\tau_1,\tau_2){\times}\Omega)\eqldef
\bigl\{ \theta\in\L^{\infty} (\tau_1,\tau_2; \W^{1,2}(\Omega)):  \ 
\dot\theta \in \L^2(\tau_1,\tau_2; \L^2(\Omega)) \bigr\},
\end{equation*}
\(0 
\le \tau_1 < \tau_2 \le T\),
endowed with the norm
\begin{equation*}
\forall \theta \in \V((\tau_1,\tau_2){\times}\Omega): \ 
\norm[\V((\tau_1,\tau_2){\times}\Omega)]{\theta}  \eqldef  
\norm[\L^{\infty}(\tau_1,\tau_2; 
\W^{1,2}(\Omega))]{\theta}
+ \norm[\L^2(\tau_1,\tau_2;\L^2(\Omega))]{\dot \theta} .
\end{equation*}
Then it follows from \eqref{eq:7} that there exists a 
generic constant $C_{\theta}>0$, independent of $\tau$,  
such that the solution of problem 
\eqref{eq:19}--\eqref{eq:21} satisfies
\begin{equation}
\label{eq:9}
\norm[\V((0,\tau){\times}\Omega)]{\theta}\leq  
C_{\theta} \exp\bigl(\tfrac{\tau}{c^c}\bigr)
\bigl(\norm[\W^{1,2}(\Omega)]{\theta^0} 
{+} \norm[\L^2(0,\tau; \L^2(\Omega))]{f^{\tilde \theta}}\bigr). 
\end{equation}
By Proposition \ref{sec:existence}, for any 
\(\tilde\theta \) belonging to a bounded subset of 
\(\C^0([0,\tau];\L^4(\Omega))\),
\(f^{\tilde \theta}\eqldef
\bfA\ee(\dot u){:}\ee(\dot u)
+\tilde\theta(\beta\bfI{:}\ee(\dot u){+}\partial_{\ze}H_2(\ze){.}\dze)
+\bfB\dze{.}\dze +\Psi(\dze)\) belongs to a  bounded subset of  
\(\L^{q/4}(0,\tau;\L^2(\Omega))\)
for any \(q>8\). Furthermore H\"older's inequality gives
\begin{equation*}
\norm[\L^2(0,\tau;\L^2(\Omega))]{f^{\tilde \theta}}
\leq \tau^{\frac{q{-}8}{2q}}\norm[\L^{q/4}(0,\tau;\L^2(\Omega))]{f^{\tilde \theta}}.
\end{equation*}
We insert this inequality into \eqref{eq:9}, we find
\begin{equation*}
\norm[\V((0,\tau){\times}\Omega)]{\theta}\leq  
C_{\theta} \exp\bigl(\tfrac{\tau}{c^c}\bigr)
\bigl(\norm[\W^{1,2}(\Omega)]{\theta^0} 
{+}\tau^{\frac{q{-}8}{2q}}\norm[\L^{q/4}(0,\tau;
\L^2(\Omega))]{f^{\tilde \theta}}\bigr). 
\end{equation*}
Thus it is clear that
\(\Phi_{\tau}^{\tilde\theta,\theta}\) maps any bounded subset of 
\(\C^0([0,\tau];\L^4(\Omega))\) into a bounded subset of
$\V((0,\tau){\times}\Omega)$. However $\V((0,\tau){\times}\Omega)$ is compactly
embedded into \(\C^0([0,\tau];\L^4(\Omega))\)  (see \cite{Sim87CSLP}), which
allows us to conclude.
\eproof

\begin{proposition}
\label{prop:Phi2} 
Let $\tau \in (0, T]$.
Assume that \eqref{eq:Psi}, \eqref{eq:H}, \eqref{eq:47}, 
\eqref{eq:LM}, \eqref{eq:Q}, \eqref{eq:f}, \eqref{eq:36}, 
\(\theta^0\in\W^{1,2}_{\kappa, \rm{Neu}}(\Omega)\), 
\(u^0\in\W^{1,4}_{\rm{Dir}}(\Omega)\)
and \(\zeo\in\W^{2,4}_{\rm{Neu}}(\Omega)\) if \(\alpha>0\) and
\(\zeo\in\L^4(\Omega)\) if \(\alpha=0\) hold.
Then the mapping
\(\Phi_{\tau}^{\tilde\theta,\theta}\) 
is continuous from  \(\C^0([0,\tau];\L^4(\Omega))\) 
into \(\C^0([0,\tau];\L^4(\Omega))\).
\end{proposition}

\bproof
Let us consider a converging sequence $(\tilde \theta_n)_{n \in \en}
\subset (\C^0([0,\tau];\L^4(\Omega)))^{\en}$ and let $\tilde 
\theta_*$ be its limit. We denote by $(u_n, \ze_n)$ the solution of 
\eqref{eq:ent_eq_1_1}--\eqref{eq:init_cond2} with 
$\tilde \theta=\tilde \theta_n$,  and 
$\theta_n\eqldef\Phi_{\tau}^{\tilde\theta,\theta}(\tilde\theta_n)$ 
for all 
$n\ge  0$. Similarly, let $(u_*, \ze_*)$ be the solution of 
\eqref{eq:ent_eq_1_1}--\eqref{eq:init_cond2} with $\tilde 
\theta=\tilde \theta_*$, and $\theta_* \eqldef\Phi_{\tau}^{\tilde\theta,\theta}
(\tilde \theta_*)$. 
Since $(\tilde \theta_n)_{n \in \en}$ is a bounded family 
of $\C^0([0,\tau];\L^4(\Omega))$, 
we infer that $(\theta_n)_{n \in \en}$ 
is bounded in $\V((0,\tau){\times}\Omega)$. 
It follows that $(\theta_n)_{n \in \en}$ is relatively compact 
in $\C^0([0,\tau];\L^4(\Omega))$ (see \cite{Sim87CSLP}). Hence, 
there exists a subsequence, still denoted by $(\theta_n)_{n \in \en}$, 
such that
\begin{equation*}
\begin{aligned}
\theta_n &\rightarrow \theta \ \text{ in } \ 
\C^0([0,\tau];\L^4(\Omega)),\\
\theta_n &\rightharpoonup \theta \ \text{ in } \ \L^{2}(0,\tau;\W^{1,2}(\Omega))
\ \text{ weak},\\
\dot\theta_n &\rightharpoonup \dot\theta \ \text{ in } \ 
\L^2(0,\tau;\L^2(\Omega))\ \text{ weak},
\end{aligned}
\end{equation*}
and for all $n \ge 0$, we have $\theta_n ( \cdot, 0)= \theta^0$ and 
\begin{equation}
\label{cont-phi-1}
\begin{aligned}
&
\int_{\calQ_{\tau}} c(x) \dot \theta_n (x,t)  \xi (x)  w  (t) \dd x \dd t 
+ \int_{\calQ_{\tau}} \kappa(x)\nabla \theta_n (x,t) {.}\nabla \xi (x) 
w(t) \dd x
\dd t \\&  
= \int_{\calQ_{\tau}} f^{\tilde \theta_n} (x,t) \xi (x) w (t) \dd x \dd t
\end{aligned}
\end{equation}
for all \(\xi\in\W^{1,2}(\Omega)\) 
and $w \in {\mathcal D}( 0, \tau)$.  We observe that  
since $(\theta_n)_{n \in \en}$ 
converges strongly to $\theta$ in $\C^0([0,\tau];\L^4(\Omega))$, 
we have immediately $\theta ( \cdot, 0) = \theta^0$.
In order to pass to the limit in \eqref{cont-phi-1}, it remains to study 
the convergence of $(f^{\tilde \theta_n} )_{n \in \en}$. We begin with the 
study of the convergence of $(u_n, \ze_n)_{n \in \en}$.

It is convenient to introduce the following functional space 
\(\X^{\alpha}(\Omega){\eqldef}\W^{1,2}_{\textrm{Neu}}(\Omega)\) if \(\alpha>0\)
and \(\X^{\alpha}(\Omega)\eqldef\L^2(\Omega)\) if \(\alpha=0\).

\begin{lemma}
\label{lm:cont_uz}
Let $\tau \in (0, T]$. Assume that \eqref{eq:Psi}, \eqref{eq:H}, 
\eqref{eq:47}, \eqref{eq:LM}, \eqref{eq:Q}, \eqref{eq:f}, 
 \(u^0\in\W^{1,4}_{\rm{Dir}}(\Omega)\)
and \(\zeo\in\W^{2,4}_{\rm{Neu}}(\Omega)\) if \(\alpha>0\) and
\(\zeo\in\L^4(\Omega)\) if \(\alpha=0\) hold.
Then the mapping \(\tilde \theta\mapsto(u,z)\), where $(u,z)$ 
is the unique solution of 
\eqref{eq:ent_eq_1_1}--\eqref{eq:init_cond2}, 
is continuous from
\(\C^0([0,\tau];\L^4(\Omega))\) into
\(\W^{1,2}(0,\tau;\W^{1,2}_{\rm{Dir}}(\Omega)\times\L^2(\Omega))
\cap\L^{\infty}(0,\tau;\W^{1,2}_{\rm{Dir}}(\Omega)\times\X^{\alpha}(\Omega))\).
\end{lemma}

\bproof
We consider \(\tilde \theta_i\in \C^0([0,\tau];\L^4(\Omega))\) and 
for $i=1,2$, we denote by \((u_i,z_i)\) the solution of the
following system:
\begin{subequations}
\label{eq:61}
\begin{align}
&
-\dive(\bfE(\ee(u_i){-}\bfQ\ze_i){+}\beta \tilde \theta_i\bfI{+}
\bfA\ee(\dot u_i))
=f,
\label{eq:62}\\&
\partial\Psi(\dze_i)+\bfB\dze_i-\tilde\bfQ^{\tra}\bfE(\ee(u_i){-}\bfQ\ze_i)
+\partial_{\ze} H_1(\ze_i)+
\tilde \theta_i \partial_{\ze}H_2(\ze_i)-\alpha\Delta\ze_i{\ni} 0,\label{eq:63}
\end{align}
\end{subequations}
with boundary conditions
\begin{equation}
\label{eq:65}
u_i=0,\quad
\alpha \nabla \ze_i{\cdot}\bfn=0 \quad\text{ on }\quad
\partial\Omega\times [0,\tau),
\end{equation}
and initial conditions
\begin{equation}
\label{eq:64}
u_i(\cdot,0)=u^0,\quad
\ze_i(\cdot,0)=\zeo \quad\text{ in }\quad \Omega.
\end{equation}
On the one hand, with the definition of the subdifferential 
\(\partial \Psi(\dze_i)\) (see \cite{Brez73OMMS}),
we have 
\begin{equation}
\label{eq:66bis}
\begin{aligned}
&
\int_{\Omega}
- \bfE(\ee(u_i){-}\bfQ\ze_i){:}( \tilde \bfQ\dze_{3{-}i} {-}
\tilde \bfQ\dze_i)\dd x + \int_{\Omega}
\bfB\dze_i{.}(\dze_{3{-}i}{-}\dze_i)\dd x\\& 
-\alpha
\int_{\Omega}
\Delta {\ze}_i{.}(\dze_{3{-}i}{-}\dze_i)\dd x
+ \int_{\Omega} 
\partial_{\ze}H_1(\ze_i){.}(\dze_{3{-i}}{-}\dze_i)\dd x\\&
+
\int_{\Omega} 
\tilde \theta_i\partial_{\ze}H_2(\ze_i){.}(\dze_{3{-i}}{-}\dze_i)\dd x
+
\int_{\Omega}  \Psi(\dze_{3{-i}}) \dd x 
-\int_{\Omega} \Psi(\dze_i) \dd x \geq 0 
\end{aligned}
\end{equation}
for almost every  \(t\in[0,\tau]\).
On the other hand, we multiply \eqref{eq:62} by 
\(\dot u_{3{-}i}{-}\dot u_i\),  we integrate this expression 
over $\Omega$ and we add it to  to \eqref{eq:66bis}. We obtain 
\begin{equation}
\label{eq:66}
\begin{aligned}
&
\int_{\Omega}
\bfE(\ee(u_i){-}\bfQ\ze_i){:}((\ee(\dot u_{3{-}i}){-}\tilde \bfQ\dze_{3{-}i}){-}
(\ee(\dot u_i){-}\tilde \bfQ\dze_i))\dd x
\\&+
\beta\int_{\Omega}
\tilde \theta_i
\bfI{:}(\ee(\dot u_{3{-}i}){-}\ee(\dot u_i))\dd x
+
\int_{\Omega}
\bfA\ee(\dot u_i){:}
(\ee(\dot u_{3{-}i}){-}\ee(\dot u_i))\dd x
\\&+
\int_{\Omega}
\bfB\dze_i{.}(\dze_{3{-}i}{-}\dze_i)\dd x  
-\alpha
\int_{\Omega}
\Delta {\ze}_i{.}(\dze_{3{-}i}{-}\dze_i)\dd x\\&
+
\int_{\Omega}
\partial_{\ze} H_1(\ze_i){.}(\dze_{3{-}i}{-}\dze_i)\dd x 
+
\int_{\Omega} 
\tilde \theta_i\partial_{\ze}H_2(\ze_i){.}(\dze_{3{-i}}{-}\dze_i)\dd x
\\&-
\int_{\Omega}
f{\cdot}(\dot u_{3{-}i}{-}\dot u_i)\dd x+
\int_{\Omega}  \Psi(\dze_{3{-i}}) \dd x - \int_{\Omega}  
\Psi(\dze_i) \dd x \geq 0 
\end{aligned}
\end{equation}
for almost every  \(t\in[0,\tau]\).
Therefore, we 
take \(i=1,2\) in \eqref{eq:66}, and thus we add these two inequalities, 
we obtain
\begin{equation}
\label{eq:11}
\begin{aligned}
&
\int_{\Omega}
\bfE(\ee(\bar u){-}\tilde \bfQ\bar\ze){:}(\ee(\dot{\bar u}){-}
\tilde \bfQ\dot{\bar\ze})\dd x
+
\int_{\Omega}
\bfA\ee(\dot{\bar u}){:}\ee(\dot{\bar u})\dd x
+
\int_{\Omega}
\bfB\dot{\bar\ze}{.}\dot{\bar\ze}\dd x\\&
- 
\alpha\int_{\Omega}\Delta {\bar\ze}{.}\dot{\bar\ze}\dd x
+
\int_{\Omega}
(\partial_{\ze} H_1(\ze_1){-}\partial_{\ze} H_1(\ze_2)){.}\dot{\bar\ze}\dd x\\&
\leq -
\beta\int_{\Omega}\bar\theta
\bfI{:}\ee(\dot{\bar u})\dd x
-
\int_{\Omega}
(\tilde \theta_1\partial_{\ze}H_2(\ze_1){-}
\tilde \theta_2\partial_{\ze}H_2(\ze_2)){.}
\dot{\bar\ze}\dd x
\end{aligned}
\end{equation}
with
\(\bar u\eqldef u_1{-}u_2\), \(\bar z\eqldef z_1{-}z_2\)
and \(\bar\theta\eqldef  \tilde \theta_1{-} \tilde \theta_2\).
Let $C^{H_1}>0$ and define
\begin{equation} 
\label{eqlp:cont_phi1}
\delta_{\alpha}(t)\eqldef
\tfrac12 \int_{\Omega}
\bfE(\ee( \bar u){-}\tilde \bfQ \bar \ze){:}
(\ee( \bar u){-} \tilde \bfQ \bar\ze)
\dd x
-  \tfrac{\alpha}2 \int_{\Omega}\Delta \bar\ze{.}\bar\ze
\dd x+ \tfrac{C^{H_1}}{2} \int_{\Omega} \abs{\bar\ze}^2 \dd x
\end{equation}
for all \(t \in [0,\tau]\). 
By using assumptions \eqref{eq:47} and \eqref{eq:Q} combined with 
Korn's and Young's inequalities, 
we find that there exists $c^{\delta}>0$ such that 
\begin{equation} 
\label{eq:14bis}
\forall t{\in}[0,\tau]:
\delta_{\alpha}(t)\ge c^{\delta}
\bigl(\norm[\W^{1,2}(\Omega)]{\bar u(\cdot,t)}^2  
{+}\norm[\L^2(\Omega)]{\bar z(\cdot,t)}^2 \bigr)  
{+} \tfrac{\alpha}{2} \norm[\L^2(\Omega)]{\nabla \bar z(\cdot,t)}^2.
\end{equation}
Furthermore Proposition \ref{sec:existence} implies that
the mapping  $\delta_{\alpha}(\cdot)$ 
is continuous on $[0,\tau]$ and its derivative 
in the sense of distributions belongs to $\L^1(0,\tau)$. Then
$\delta_{\alpha}(\cdot)$ 
is absolutely continuous on $[0,\tau]$ and with \eqref{eq:H3}, 
\eqref{eq:LM}, \eqref{eq:11} and \eqref{eqlp:cont_phi1}, we get
\begin{equation}
\label{eq:1}
\begin{aligned}
&
\dot\delta_{\alpha}(t)
+c^{\bfA}
\norm[\L^2(\Omega)]{\ee(\dot{\bar u})}^2 
+c^{\bfB}
\norm[\L^2(\Omega)]{\dot{\bar\ze}}^2
\le 
(C_{\ze \ze}^{H_1} {+}  C^{H_1})  
\int_{\Omega} \abs{\bar \ze} \abs{\dot{\bar \ze}} \dd x\\& 
-\beta\int_{\Omega}\bar\theta
\bfI{:}\ee(\dot{\bar u})\dd x
-\int_{\Omega}(\tilde \theta_1\partial_{\ze}H_2(\ze_1){-}
\tilde \theta_2\partial_{\ze}H_2(\ze_2)){.}\dot{\bar \ze}\dd x
\end{aligned}
\end{equation}
for almost every \(t\in [0,\tau]\).

Let us distinguish now the cases $\alpha=0$ and $\alpha >0$. 

If $\alpha =0$, then $\partial_{\ze} H_2 \equiv 0$ and \eqref{eq:1} reduces to 
\begin{equation*}
\dot\delta_{0}(t)
+c^{\bfA}
\norm[\L^2(\Omega)]{\ee(\dot{\bar u})}^2 
+c^{\bfB}
\norm[\L^2(\Omega)]{\dot{\bar\ze}}^2
\le 
(C_{\ze \ze}^{H_1}{+}C^{H_1})  \int_{\Omega}
 \abs{\bar \ze} \abs{\dot{\bar \ze}} \dd x
-\beta\int_{\Omega}\bar\theta
\bfI{:}\ee(\dot{\bar u})\dd x
\end{equation*}
for almost every \(t\in [0,\tau]\).
The two  terms of the right hand side can be estimated 
by using Cauchy-Schwarz's inequality, we have
\begin{equation*}
\begin{aligned}
\dot\delta_{0}(t)
+\tfrac{c^{\bfA}}2
\norm[\L^2(\Omega)]{\ee(\dot{\bar u})}^2 
+\tfrac{c^{\bfB}}2
\norm[\L^2(\Omega)]{\dot{\bar\ze}}^2
&\leq
\tfrac{9 \beta^2}{2c^{\bfA}}\norm[\L^2(\Omega)]{\bar\theta}^2+
\tfrac{(C_{\ze \ze}^{H_1}  +  C^{H_1})^2}{2c^{\bfB}}\norm[\L^2(\Omega)]{\bar \ze}^2
\\&\leq \tfrac{9 \beta^2}{2c^{\bfA}}\norm[\L^2(\Omega)]{\bar\theta}^2+
\tfrac{(C_{\ze \ze}^{H_1}  + C^{H_1})^2}{2c^{\bfB}c^{\delta }}\delta_{\alpha}(t).
\end{aligned}
\end{equation*}
Therefore we integrate over \((0,t)\) and we use Gr\"onwall's
lemma, we find
\begin{equation*}
\begin{aligned}
&\delta_{0}(t)
+\tfrac{c^{\bfA}}2
\norm[\L^2(0, t; \L^2(\Omega))]{\ee(\dot{\bar u})}^2 
+\tfrac{c^{\bfB}}2
\norm[\L^2(0, t ; \L^2(\Omega))]{\dot{\bar\ze}}^2
\\&\leq
\tfrac{9 \beta^2}{2c^{\bfA}} t \| \bar\theta \|^2_{\C^0([0, \tau]; \L^2(\Omega))}  
\exp\bigl(\tfrac{(C_{\ze \ze}^{H_1}  + C^{H_1})^2}{2c^{\bfB}c^{\delta}} \tau
\bigr)
\end{aligned}
\end{equation*}
for all $t \in [0, \tau]$.

If $\alpha \not = 0$,
the following decomposition is used to
estimate the last term in \eqref{eq:1}, namely  
\begin{equation*}
(\tilde \theta_1\partial_{\ze}
H_2(\ze_1){-} \tilde \theta_2\partial_{\ze}H_2(\ze_2)){.}
\dot{\bar\ze}=(\bar\theta\partial_{\ze}H_2(\ze_1){+}
\tilde \theta_2(\partial_{\ze}H_2(\ze_1){-}\partial_{\ze}H_2(\ze_2))){.}
\dot{\bar \ze}.
\end{equation*}
Then it follows that
\begin{equation}
\label{eq:71}
\begin{aligned}
&\dot \delta_{\alpha}(t) 
+\tfrac{c^{\bfA}}2
\norm[\L^2(\Omega)]{\ee(\dot{\bar u})}^2 
+\tfrac{3c^{\bfB}}4
\norm[\L^2(\Omega)]{\dot{\bar\ze}}^2\\&
\leq
\tfrac{9 \beta^2}{2c^{\bfA}}
\norm[\L^2(\Omega)]{\bar \theta}^2+
 \tfrac{(C_{\ze \ze}^{H_1} +  C^{H_1})^2}{c^{\bfB}}
\norm[\L^2(\Omega)]{\bar\ze}^2
\\&+\int_{\Omega}
\big(\abs{\bar\theta}\abs{\partial_{\ze}H_2(\ze_1)}\abs{\dot{\bar\ze}}
{+}\abs{\tilde \theta_2}\abs{\partial_{\ze}
H_2(\ze_1){-}\partial_{\ze}H_2(\ze_2)}\abs{\dot{\bar\ze}}\bigr)\dd x.
\end{aligned}
\end{equation}
Observe that \eqref{eq:H3}, \eqref{eq:H2} and Young's inequality give 
\begin{equation*}
\begin{aligned}
&\int_{\Omega} 
\big(\abs{\bar\theta}\abs{\partial_{\ze}H_2(\ze_1)}
\abs{\dot{\bar\ze}}{+}\abs{\tilde \theta_2}\abs{\partial_{\ze}
H_2(\ze_1){-}\partial_{\ze}H_2(\ze_2)}\abs{\dot{\bar\ze}}\bigr)\dd
x
\\&\leq 
C_{\ze}^{H_2}
\int_{\Omega}
\bigl(1{+}\abs{\ze_1}\bigr)\abs{\bar\theta}\abs{\dot{\bar\ze}}
\dd x
+C_{\ze\ze}^{H_2}
\int_{\Omega} 
\abs{\tilde \theta_2}\abs{\bar\ze}\abs{\dot{\bar\ze}}\dd x
\leq
\tfrac{C_{\ze}^{H_2}}{2\gamma_1}
\norm[\L^2(\Omega)]{\bar\theta}^2 \\&+
\tfrac{C_{\ze}^{H_2}}{2\gamma_2}
\int_{\Omega}
\abs{\bar\theta}^2\abs{\ze_1}^2\dd x +
\tfrac{C_{\ze\ze}^{H_2}}{2\gamma_3}
\int_{\Omega}
\abs{\tilde \theta_2}^2\abs{\bar\ze}^2\dd x
+\tfrac{C_{\ze}^{H_2}(\gamma_1{+}\gamma_2){+}C_{\ze\ze}^{H_2}\gamma_3}2
\norm[\L^2(\Omega)]{\dot{\bar z}}^2,
\end{aligned}
\end{equation*}
with \(\gamma_i>0\),
\(i=1,2,3\).
Note that $\ze_1 \in \L^{q/2} (0,\tau; \W^{2,4}_{\textrm{Neu}} 
(\Omega))$ and 
$\W^{2,4}_{\textrm{Neu}} (\Omega) 
\hookrightarrow \L^{\infty}(\Omega)$ with continuous embedding, thus
by H\"older's inequality we have
\begin{equation}
\label{eq:73}
\begin{aligned}
&\int_{\Omega} 
\big(\abs{\bar\theta}\abs{\partial_{\ze}H_2(\ze_1)}
\abs{\dot{\bar\ze}}{+}\abs{\theta_2}\abs{\partial_{\ze}
H_2(\ze_1){-}\partial_{\ze}H_2(\ze_2)}\abs{\dot{\bar \ze}}\bigr)\dd x
\\&\leq \tfrac{C_{\ze}^{H_2}}{2\gamma_1}
\norm[\L^2(\Omega)]{\bar\theta}^2 
+
\tfrac{C_{\ze}^{H_2}}{2 \gamma_2}
\norm[\L^{\infty}(\Omega)]{\ze_1}^2
\norm[\L^2(\Omega)]{\bar\theta}^2\\&+
\tfrac{C_{\ze\ze}^{H_2}}{2\gamma_3}
\norm[\L^4(\Omega)]{\tilde \theta_2}^2
\norm[\L^4(\Omega)]{\bar \ze}^2
{+}\tfrac{ C_{\ze}^{H_2}(\gamma_1{+}\gamma_2){+}C_{\ze\ze}^{H_2}\gamma_3}2
\norm[\L^2(\Omega)]{\dot{\bar \ze}}^2.
\end{aligned}
\end{equation}
We insert \eqref{eq:73} in \eqref{eq:71} and we choose
\(\gamma_1=\gamma_2=\tfrac{c^{\bfB}}{4C_{\ze}^{H_2}}\) and 
\(\gamma_3=\tfrac{c^{\bfB}}{2C_{\ze\ze}^{H_2}}\). Using
the continuous embedding 
\(\W^{1,2}(\Omega)\hookrightarrow\L^4(\Omega)\) and \eqref{eq:14bis}, we obtain 
\begin{equation}
\label{eq:74}
\begin{aligned}
&\dot \delta_{\alpha} (t) 
{+}\tfrac{c^{\bfA}}2
\norm[\L^2(\Omega)]{\ee(\dot{\bar u})}^2 
{+}\tfrac{c^{\bfB}}4
\norm[\L^2(\Omega)]{\dot{\bar \ze}}^2\\&
\leq \bigl(\tfrac{9\beta^2}{2c^{\bfA}}{+}
\tfrac{2(C_{\ze}^{H_2})^2}{c^{\bfB}}{+}
 \tfrac{2 (C_{\ze}^{H_2})^2}{c^{\bfB}}
\norm[\L^{\infty}(\Omega)]{\ze_1}^2
\bigr)
\norm[\L^2(\Omega)]{\bar \theta}^2 \\&+
\tfrac1{c^{\delta}_{\alpha}c^{\bfB}}
\bigl((C_{\ze \ze}^{H_1}{+}C^{H_1})^2{+}C_1(C_{\ze\ze}^{H_2})^2
\norm[\L^4(\Omega)]{\tilde \theta_2}^2\bigr)
\delta_{\alpha} (t)
\end{aligned}
\end{equation}
for almost every \(t \in [0,\tau]\), where $C_1$ is the generic 
constant involved in the continuous embedding of $\W^{1,2}(\Omega)$ 
into $\L^4(\Omega)$ and $c^{\delta}_{\alpha}\eqldef 
\min\bigl(c^{\delta}, \tfrac{\alpha}{2}\bigr)$. Let us define
\begin{equation*}
c (\tilde \theta_2)\eqldef
\tfrac1{c^{\delta}_{\alpha} c^{\bfB}}
\bigl((C_{\ze \ze} ^{H_1}{+} C^{H_1})^2{+}C_1 (C_{\ze\ze}^{H_2})^2
\| \tilde \theta_2 \|^2_{\C^0([0, \tau], \L^4(\Omega))} \bigr).
\end{equation*}
By using Gr\"onwall's lemma, we get 
\begin{equation*}
\begin{aligned}
& \delta_{\alpha}(t) +\tfrac{c^{\bfA}}2
\norm[\L^2(0,t; \L^2(\Omega))]{\ee(\dot{\bar u})}^2 
+\tfrac{c^{\bfB}}4
\norm[\L^2(0, t; \L^2(\Omega))]{\dot{\bar \ze}}^2 \\
&
\leq  \bigl(\tfrac{9 \beta^2}{2c^{\bfA}} \tau {+}
\tfrac{2(C_{\ze}^{H_2})^2}{c^{\bfB}}  \tau {+}
\tfrac{2 (C_{\ze}^{H_2})^2}{c^{\bfB}}
\norm[\L^2(0, \tau; \L^{\infty}(\Omega))]{\ze_1}^2 \bigr)
\| \bar \theta \|^2_{\C^0([0,\tau];\L^2(\Omega))}\\&
(1{+}\tau c(\tilde\theta_2) \exp ( c (\tilde \theta_2) \tau) )
\end{aligned}
\end{equation*}
for all \(t\in[0,\tau]\). 
\eproof

As a corollary, it is possible to prove that

\begin{lemma}
\label{lm:cont_ftheta}
Let $\tau \in (0, T]$. Assume that \eqref{eq:Psi}, 
\eqref{eq:H}, \eqref{eq:47}, \eqref{eq:LM}, \eqref{eq:Q}, \eqref{eq:f}, 
 \(u^0\in\W^{1,4}_{\rm{Dir}}(\Omega)\)
and \(\zeo\in\W^{2,4}_{\rm{Neu}}(\Omega)\) if \(\alpha>0\) and
\(\zeo\in\L^4(\Omega)\) if \(\alpha=0\) hold.
Then the mapping \(\tilde \theta \mapsto f^{\tilde \theta}\) with 
\( f^{\tilde \theta}= \bfA\ee(\dot u){:}\ee(\dot u)
+\tilde\theta(\beta\bfI{:}\ee(\dot u){+}\partial_{\ze}H_2(\ze){.}\dze)
+\bfB\dze{.}\dze +\Psi(\dze)
\), where $(u,z)$ is the unique solution of 
\eqref{eq:ent_eq_1_1}--\eqref{eq:init_cond2}, 
is continuous from
\(\C^0([0,\tau];\L^4(\Omega))\) into $\L^{r}(0,\tau; \L^{4/3}( \Omega))$, with 
$\frac{1}{ r}= \frac{2}{q} + \frac12$.
\end{lemma}

\bproof
We consider once again \(\tilde \theta_i\in \C^0([0,\tau];\L^4(\Omega))\) and 
for $i=1,2$, we denote by \((u_i,z_i)\) the solution of the
 system \eqref{eq:61}--\eqref{eq:64}.
With the definition of $f^{\tilde \theta}$ we have
\begin{equation*}
\begin{aligned}
&f^{\tilde \theta_1}-f^{\tilde \theta_2}= 
\bfA\ee(\dot u_1{+}\dot u_2){:}\ee(\dot u_1{-}\dot u_2)+
(\tilde\theta_1{-}\tilde\theta_2)(\beta\bfI{:}\ee(\dot u_1){+}
\partial_{\ze} H_2(\ze_1){.}\dot \ze_1)
\\&
+\tilde\theta_2(\beta 
\bfI{:}\ee(\dot u_1{-}\dot u_2){+}\partial_{\ze}H_2(\ze_1){.}\dot \ze_1{-}
\partial_{\ze}H_2(\ze_2){.}
\dot \ze_2)\\& 
{+}\bfB (\dot \ze_1{+}\dot \ze_2){.}
(\dot \ze_1{-}\dot \ze_2)
{+}\Psi(\dot \ze_1){-}\Psi(\dot \ze_2).
\end{aligned}
\end{equation*}
Thus we have
\begin{equation*}
\begin{aligned}
&\abs{f^{\tilde \theta_1}{-}f^{\tilde \theta_2}} \le 
\norm{\bfA} \abs{\ee(\dot u_1{+}\dot u_2)} \abs{\ee(\dot{\bar u})} +
\abs{\bar \theta} (3 \beta \abs{\ee(\dot u_1)} {+} C_z^{H_2} (1{+}\abs{\ze_1}) 
\abs{\dot \ze_1})
\\&
+\abs{\tilde\theta_2} (3 \beta 
\abs{\ee(\dot{\bar u})} {+}\abs{\partial_{\ze}H_2(\ze_1){.}\dot \ze_1{-}
\partial_{\ze}H_2(\ze_2){.}
\dot \ze_2} ) \\&
+ \norm{\bfB} \abs{\dot \ze_1{+}\dot \ze_2}
\abs{\dot{\bar \ze}}
+\abs{\Psi(\dot \ze_1){-}\Psi(\dot \ze_2)}.
\end{aligned}
\end{equation*}
But \eqref{eq:Psi.bdd} and \eqref{eq:4} lead to
\begin{equation*}
\abs{\Psi(\dot \ze_1){-}\Psi(\dot \ze_2)} \le C^{\Psi} 
\abs{\dot \ze_1{-}\dot \ze_2} = C^{\Psi} \abs{\dot{\bar \ze}}
\end{equation*}
and  \eqref{eq:H2} and \eqref{eq:H3} give
\begin{equation*}
\begin{aligned}
\abs{\partial_{\ze}H_2(\ze_1){.}\dot \ze_1{-}\partial_{\ze}
H_2(\ze_2){.}\dot \ze_2 }
&\leq\abs{\partial_{\ze}H_2(\ze_1)}\abs{\dot{\bar  \ze}}+
\abs{\partial_{\ze}H_2(\ze_1){-}\partial_{\ze}H_2(\ze_2) } 
\abs{\dot \ze_2}
\\&
\leq C_{\ze}^{H_2} (1{+}\abs{\ze_1}) 
\abs{\dot{\bar  \ze}}+C^{H_2}_{\ze\ze}\abs{\bar \ze}
\abs{\dot \ze_2}. 
\end{aligned}
\end{equation*}
Then, recalling that $\partial_{\ze} H_2 \equiv 0$ whenever $\alpha =0$, 
the boundedness properties of $(u,z)$ stated at Proposition 
\ref{sec:existence} and the continuity property of $\tilde \theta 
\mapsto (u,z)$ proved in  Lemma \ref{lm:cont_uz}  allow us to conclude 
by using Young's inequality.
\eproof

Now we may conclude the proof of Proposition \ref{prop:Phi2}. Indeed, 
since $(\tilde \theta_n)_{n \in \en}$ converges strongly to 
$\tilde\theta_*$ in $\C^0([0,\tau];\L^4(\Omega))$, 
we infer from Lemma \ref{lm:cont_ftheta} that 
\begin{equation*}
\lim_{n \rightarrow + \infty}\int_{\calQ_{\tau}}f^{\tilde \theta_n}(x,t) 
\xi(x)w(t)\dd x\dd t=\int_{\calQ_{\tau}}f^{\tilde\theta_*}(x,t) 
\xi(x)w (t)\dd x\dd t
\end{equation*}
for all \(\xi\in\W^{1,2}(\Omega)\) 
and \(w \in {\mathcal D} (0, \tau) \).
Therefore we may pass to the limit in all the terms  of \eqref{cont-phi-1} 
to get
\begin{equation}
\label{eq:10}
\begin{aligned}
&
\int_{\calQ_{\tau}} c(x) \dot \theta (x,t)  \xi (x)  w  (t) \dd x \dd t 
+ \int_{\calQ_{\tau}} \kappa(x) \nabla \theta (x,t) \nabla \xi (x) 
w(t) \dd x
\dd t \\
& = 
\int_{\calQ_{\tau}} f^{\tilde \theta_*} (x,t)\xi(x)w(t)\dd x\dd t
\end{aligned}
\end{equation}
for all \(\xi\in\W^{1,2}(\Omega)\) 
and \(w \in {\mathcal D} (0, \tau) \).
It follows  that $\theta$ is solution of problem 
\eqref{eq:19}--\eqref{eq:21} 
with the data $f^{\tilde \theta_*}$. Besides by uniqueness of 
the solution, we have $\theta=\theta_*$ and the 
whole sequence $(\theta_n)_{n \in \en}$ converges to 
$\theta_*$ in $\C^0([0,\tau];\L^4(\Omega))$.
\eproof

\begin{corollary} 
\label{corollaire}
Let $\tau \in (0, T]$. Assume that \eqref{eq:Psi}, \eqref{eq:H}, 
\eqref{eq:47}, \eqref{eq:LM}, \eqref{eq:Q}, \eqref{eq:f},  
\eqref{eq:36}, \(\theta^0\in\W^{1,2}_{\kappa,\rm{Neu}}(\Omega)\), 
 \(u^0\in\W^{1,4}_{\rm{Dir}}(\Omega)\)
and \(\zeo\in\W^{2,4}_{\rm{Neu}}(\Omega)\) if \(\alpha>0\) and
\(\zeo\in\L^4(\Omega)\) if \(\alpha=0\) hold.
Then there exists $\tau \in (0,T]$ such that 
\(\Phi_{\tau}^{\tilde\theta,\theta}\) admits a fixed point in 
$\C^0([0,\tau];\L^4(\Omega))$.
\end{corollary}

\bproof
We have already proved in the previous propositions that 
\(\Phi_{\tau}^{\tilde\theta,\theta}\) is a continuous mapping from 
$\C^0([0,\tau];\L^4(\Omega))$ into $\C^0([0,\tau];\L^4(\Omega))$ 
and maps any bounded subset ${\mathcal C} \subset \C^0([0,\tau];\L^4(\Omega))$ 
into a  bounded relatively compact subset. Hence we will be able to conclude 
by using Schauder's fixed point theorem (see \cite{Evan10PDES}) 
if we can find a closed convex bounded subset ${\mathcal C}$ of
$\C^0([0,\tau];
\L^4(\Omega))$ such that $\Phi_{\tau}^{\tilde\theta,\theta} ({\mathcal C}) 
\subset {\mathcal C}$. 

Let $C_1>0$ be the generic constant involved in the continuous embedding of 
$\W^{1,2}(\Omega)$ into $\L^4(\Omega)$ and let 
$R^{\theta} > C_1 C_{\theta} \exp( \tfrac{T}{c^c}) 
\norm[\W^{1,2}(\Omega)]{\theta^0}$, where $C_{\theta}$ is 
the constant defined in Proposition \ref{prop:Phi1}. For any $\tau \in (0,T]$ 
and $\tilde \theta  \in {\mathcal C} 
\eqldef {\bar B}_{\C^0([0,\tau];\L^4(\Omega))} (0, R^{\theta})$, 
we denote $\theta = \Phi_{\tau}^{\tilde\theta,\theta} (\tilde \theta)$ 
and we have (see \eqref{eq:9})
\begin{equation*}
\norm[\V((0,\tau){\times}\Omega)]{\theta}\leq  
C_{\theta} \exp\bigl(\tfrac{\tau}{c^c}\bigr)
\bigl(\norm[\W^{1,2}(\Omega)]{\theta^0} 
{+} \norm[\L^2(0,\tau; \L^2(\Omega))]{f^{\tilde \theta}}\bigr)
\end{equation*}
and $\theta \in \C^0([0,\tau];\L^4(\Omega))$. Thus we have
\begin{equation} 
\label{eq:exist_1}
\begin{aligned}
&\| \theta\|_{\C^0([0,\tau];\L^4(\Omega))} = \norm[\L^{\infty}(0, \tau; 
\L^4(\Omega))]{\theta} 
\\&\leq 
C_1 C_{\theta} \exp\bigl(\tfrac{\tau}{c^c}\bigr)
\bigl(\norm[\W^{1,2}(\Omega)]{\theta^0} 
{+} \tau^{\frac{q-8}{2q}}\norm[\L^{q/4}(0,\tau; 
\L^2(\Omega))]{f^{\tilde \theta}}\bigr)
\end{aligned}
\end{equation}
for any $q>8$. 
Since $\lim_{\tau \to 0} \tau^{\frac{q-8}{2q}} =0$, we only need now to establish
that $\norm[\L^{q/4}(0,\tau; \L^2(\Omega))]{f^{\tilde \theta}}$ remains 
bounded independently of $\tau$. Let us emphasize that Proposition 
\ref{sec:existence} implies that $f^{\tilde \theta}$ remains in a 
bounded subset of $\L^{q/4}(0,\tau; \L^2(\Omega))$ but this does 
not allow us to conclude since we don't know if the diameter of this 
bounded subset depends on $\tau$ or not. In order to cope with this 
difficulty, we consider the extension of $\tilde \theta$ to $[0,T]$ by 
zero on $(\tau, T]$. We denote by $\tilde \theta_{\textrm{ext}}$ this 
extension. Of course, for any $\tilde \theta \in {\mathcal C}$, we 
have $\tilde \theta_{\textrm{ext}} \in \L^{q}(0, T; \L^4(\Omega))$ for 
any $q>8$ and
\begin{equation*}
\norm[\L^{q}(0,T; \L^4(\Omega))]{\tilde \theta_{\textrm{ext}}  } 
= \norm[\L^{q}(0,\tau; \L^4(\Omega))]{\tilde \theta} = \tau^{\frac1q} \| 
\tilde \theta\|_{\C^0([0,\tau];\L^4(\Omega))}  \le T^{\frac1q} R^{\theta}.
\end{equation*}
Then we define $(u_{\textrm{ext}}, \ze_{\textrm{ext}})$ as the unique 
solution of problem \eqref{eq:ent_eq_1_1}--\eqref{eq:init_cond2} with 
$\tau$ replaced by $T$ and $\tilde \theta$ replaced by  
$ \tilde \theta_{\textrm{ext}}$. Since $\tilde \theta_{\textrm{ext}}$ 
remains in the closed ball ${\bar B}_{\L^{q}(0,T; \L^4(\Omega))} (0, 
T^{1/q} R^{\theta})$, it is clear that 
Proposition \ref{sec:existence} implies that 
$(u_{\textrm{ext}}, \ze_{\textrm{ext}})$ remains in a bounded subset 
of $ \W^{1,q}(0,\tau;\W^{1,4}_{\textrm{Dir}}(\Omega)) \times  
( \L^{q/2} (0,\tau;\W^{2,4}_{\textrm{Neu}}(\Omega))\ \cap\  \C^0([0, \tau]; 
\W^{1,2}_{\textrm{Neu}}(\Omega))
\ \cap \ \W^{1,q/2}(0,\tau;\L^{4}(\Omega))
\cap\W^{1,q}(0,\tau;\L^{2}(\Omega)))$ if \(\alpha>0\) or 
$\W^{1,q}(0,\tau;\W^{1,4}_{\textrm{Dir}}(\Omega)) \times 
\W^{1,q}(0,\tau;\L^4(\Omega))$ if \(\alpha=0\). It follows that 
\begin{equation*}
f^{\tilde \theta_{\textrm{ext}}}\eqldef \bfA\ee(\dot u_{\textrm{ext}}){:}
\ee(\dot u_{\textrm{ext}})
+\tilde\theta_{\textrm{ext}}(\beta\bfI{:}\ee(\dot u_{\textrm{ext}}){+}
\partial_{\ze}H_2(\ze_{\textrm{ext}}){.}\dot\ze_{\textrm{ext}}) 
+\bfB\dot\ze_{\textrm{ext}}{.}\dot \ze_{\textrm{ext}}+\Psi(\dot\ze_{\textrm{ext}})
\end{equation*}
remains in a bounded subset of $\L^{q/4}(0,T; \L^2(\Omega))$, i.e. 
there exists a constant $C(R^{\theta})$, depending only on $R^{\theta}$ 
and the data, such that 
$\norm[\L^{q/4}(0,T; \L^2(\Omega))]{f^{\tilde \theta_{\textrm{ext}}}} 
\le C(R^{\theta})$. But $f^{\tilde \theta}$ coincide with 
$f^{\tilde \theta_{\textrm{ext}}}$ on $[0, \tau]$ and 
\begin{equation} 
\begin{aligned}
\label{eq:exist_2}
&\norm[\L^{q/4}(0,\tau; \L^2(\Omega))]{f^{\tilde \theta} }
= \norm[\L^{q/4}(0,T; \L^2(\Omega))]{f^{\tilde 
\theta_{\textrm{ext}}} {\bf 1}_{[0, \tau]} }\\&\le 
\norm[\L^{q/4}(0,T; \L^2(\Omega))]{f^{\tilde \theta_{\textrm{ext}}}}  
\le C(R^{\theta}).
\end{aligned}
\end{equation}
Then by introducing \eqref{eq:exist_2} into \eqref{eq:exist_1} and 
by choosing $\tau \in (0, T]$ such that
\begin{equation*}
C_1 C_{\theta} \exp\bigl(\tfrac{\tau}{c^c}\bigr)
\bigl(\norm[\W^{1,2}(\Omega)]{\theta^0} 
{+} \tau^{\frac{q-8}{2q}} C(R^{\theta}) ) \le R^{\theta},
\end{equation*}
we may conclude.
\eproof

We can consider $\tau \in (0,T]$ such that 
$\Phi_{\tau}^{\tilde\theta,\theta}$ admits a fixed point $\theta$ 
in the space $\C^0([0,\tau];\L^4(\Omega))$ and we define $(u, \ze)$ as 
the unique solution of problem \eqref{eq:ent_eq_1_1}--\eqref{eq:init_cond2} 
with $\tilde \theta = \theta$. By definition of 
$\Phi_{\tau}^{\tilde\theta,\theta}$, $(u,z, \theta)$ is a solution of the 
coupled problem \eqref{eq:ent_eq}--\eqref{eq:init_cond} and by combining 
the regularity results for $(u, \ze)$ given at Proposition \ref{sec:existence} 
with the regularity results for the heat-transfer equation recalled in the 
proof of Proposition \ref{prop:Phi1}, we get  
\(\theta\in\L^{\infty}(0,\tau;\W^{1,2}_{\textrm{$\kappa$, Neu}}(\Omega))
\cap\C^0(0,\tau;\L^4(\Omega))\), 
 \(\dot\theta\in\L^2(0,\tau;\L^2(\Omega))\) and
\(u\in\W^{1,q}(0,\tau;\W^{1,4}_{\textrm{Dir}}(\Omega))\), 
\(\ze\in\L^{q/2} (0,\tau;\W^{2,4}_{\textrm{Neu}}(\Omega))
\cap \C^0([0, \tau]; \W^{1,2}_{\textrm{Neu}}(\Omega))
\cap\W^{1,q/2}(0,\tau;\L^{4}(\Omega))
\cap\W^{1,q}(0,\tau;\L^{2}(\Omega))\) when \(\alpha>0\),
\(\ze\in\W^{1,q}(0,\tau;\L^4(\Omega))\)
when \(\alpha=0\), for any $q>8$. Hence the proof of Theorem 
\ref{thm:local_existence} is complete.

%%%%%%%%%%%%%%%%%%%%%%%%%%%%%%%%%%%%%%%%%%%%%%%%%%%%%%%%%%%%%%%%%%%%%%%%%%%%%%%%

\section{Further properties of the solution}
\label{sec:furth-prop-solut}

Let us recall that system \eqref{eq:ent_eq}--\eqref{eq:init_cond}  
is thermodynamically consistent if the temperature remains positive 
(see Section \ref{sec:description-problem}). So we begin this section 
by proving that the solutions $(u, \ze, \theta)$ of 
\eqref{eq:ent_eq}--\eqref{eq:init_cond}  are physically admissible, 
i.e. $\theta (x,t) >0$ almost everywhere in $\calQ_{\tau}$. 
To this aim we introduce the following assumption for the initial temperature:

%\noindent (A-8)  
\begin{enumerate}[({A}--8)]
\itemsep0.1em

\item There exists $\bar \theta>0$ such that
\begin{equation}
\label{init_temp}
\theta^0(x)\ge\bar\theta>0 \ \text{ a.e. }\ x\in\Omega.
\end{equation}
\end{enumerate}

\vspace{0.1em}
\begin{proposition}
\label{temp_adm}
Assume that \eqref{eq:Psi}, \eqref{eq:H},  
\eqref{eq:47}, \eqref{eq:LM},  
\eqref{eq:Q}, \eqref{eq:f}, \eqref{eq:36},
\(\theta^0\in\W^{1,2}_{\kappa, \rm{Neu}}(\Omega)\), 
\(u^0\in\W^{1,4}_{\rm{Dir}}(\Omega)\)
and \(\zeo\in\W^{2,4}_{\rm{Neu}}(\Omega)\) if \(\alpha>0\) and
\(\zeo\in\L^4(\Omega)\) if \(\alpha=0\) hold. Assume also that 
condition \eqref{init_temp} is satisfied and  
$\kappa \in \C^1({\bar \Omega})$. Then,  
any solution $(u, \ze, \theta) $ of problem 
\eqref{eq:ent_eq}--\eqref{eq:init_cond} 
defined on $[0,\tau]$, $\tau \in (0,T]$,  is 
thermodynamically admissible, i.e. $\theta(x,t) >0$ 
for almost every $(x,t) \in \calQ_{\tau}$.
\end{proposition}

\bproof
The key tool of the proof is the classical Stampacchia's 
truncation method (see \cite{Brez83AFTA}). 
So we consider 
a function \(\calG\in\C^1(\Er; \Er)\) 
such that
\begin{enumerate}[(i)]
\item \(\exists C^{\calG'}>0, \ \forall \sigma\in\Er:\ \abs{\calG'(\sigma)}
\leq C^{\calG'}\),
\item \(\calG\) is strictly increasing on \((0,\infty)\),
\item \(\forall \sigma\leq 0:\ \calG(\sigma)=0\),
\end{enumerate}
and we define \(\Gamma(\sigma)\eqldef\int_0^{\sigma} \calG(s)\dd s\) 
for all $\sigma\in\Er$.
Now let $(u,z, \theta)$ be a solution of 
\eqref{eq:ent_eq}--\eqref{eq:init_cond} on $[0,\tau]$. 
We will prove that $\theta$ is positive almost everywhere in $\calQ_{\tau}$ 
in two steps: first we will establish that $\theta$ is non negative, then 
that $\theta$ remains bounded from below by a positive quantity.

Since we assumed that $\kappa \in \C^1( {\bar \Omega})$, we 
can infer that $\theta \in \L^2(0, \tau, \W^{2,2}(\Omega))$. Indeed, 
$\theta$ is a fixed point of $\Phi_{\tau}^{\tilde\theta,\theta}$, thus
\begin{equation}
\label{eq:12}
- \dive (\kappa (x) \nabla\theta) = f^{\theta} - c (x) \dot \theta
\end{equation}
with $f^{\theta}= \bfA\ee(\dot u){:}\ee(\dot u)
+\theta(\beta\bfI{:}\ee(\dot u){+}\partial_{\ze}H_2(\ze){.}\dze)
+ \bfB\dze{.}\dze +\Psi(\dze) \in \L^2(0, \tau; \L^2 (\Omega))$ and 
$c (x) \dot \theta \in  \L^2(0, \tau; \L^2 (\Omega))$. It follows 
that $- \dive (\kappa (x) \nabla\theta)  \in \L^2(0, \tau; \L^2 (\Omega))$. 
We can consider the time variable as a parameter and the linearity of 
the operator $- \dive (\kappa (x) \nabla \cdot)$ combined with classical 
regularity properties (see \cite{Brez83AFTA}) yield the announced result. 

Then we introduce the mapping $\varphi: [0,\tau]\rightarrow\Er$ given by 
\begin{equation}
\label{eq:16}
\varphi(t)\eqldef\exp\Bigl(- \tfrac1{c^c} \int_0^t 
 \tfrac{9 \beta^2}{2 c^{\bfA}} \norm[\L^{\infty}
(\Omega)]{\theta(\cdot,s)}\dd s\Bigr)
\end{equation}
for all \(t\in[ 0, \tau]\) if $\alpha =0$ and by 
\begin{equation}
\label{eq:16bis}
\varphi(t)\eqldef\exp\Bigl(- \tfrac1{c^c} \int_0^t 
\bigl( \tfrac{9 \beta^2}{2 c^{\bfA}}{+}\tfrac{(C_z^{H_2})^2}{c^{\bfB}} 
\bigl(1{+}\norm[\L^{\infty}
(\Omega)]{z(\cdot,s)}^2\bigr)\bigr)\norm[\L^{\infty}
(\Omega)]{\theta(\cdot,s)}\dd s\Bigr)
\end{equation}
for all \(t\in[ 0, \tau]\)
if $\alpha >0$. Recalling that 
$z \in \L^{q/2} (0, \tau; \W^{2,4}_{\textrm{Neu}} (\Omega))$ 
for any $q>8$ if $\alpha >0$ and $\W^{2,2}(\Omega) 
\hookrightarrow \L^{\infty}( \Omega)$, we can deduce 
that $\varphi \in \W^{1,1}(0, \tau)$ and $0 \le \varphi (t) \le 1$ 
for almost every $t \in [0, \tau]$.
Next we define
\(\Theta_{\theta_{\varphi}}(t)
\eqldef\int_{\Omega}c(x)\Gamma(\theta_{\varphi}(x,t))\dd x\) with
\(\theta_{\varphi}(x,t)\eqldef -\theta(x,t)\varphi(t)\) for almost 
every $(x,t) \in \calQ_{\tau}$.
Since $\theta \in \V((0, \tau) {\times} \Omega) \cap 
\L^2 (0, \tau, \W^{2,2} (\Omega))$, we get $\theta_{\varphi} 
\in \L^{\infty} (0, \tau; \W^{1,2} (\Omega)) \cap 
\W^{1,1} (0, \tau; \L^2(\Omega))$ and 
\begin{equation*}
\begin{aligned}
& \dot\theta_{\varphi}(x,t)=
\bigl(-\dot\theta (x,t) {+}\tfrac{\theta(x,t)}{c^c}
\tfrac{9 \beta^2}{2 c^{\bfA}} \norm[\L^{\infty}
(\Omega)]{\theta(\cdot,t)} \bigr)\varphi(t)\textrm{ if } 
\alpha =0 , \\
& \dot\theta_{\varphi}(x,t)=
\bigl(-\dot\theta(x,t){+}\tfrac{\theta(x,t)}{c^c}
\bigl(\tfrac{9 \beta^2}{2 c^{\bfA}}{+}
\tfrac{(C_z^{H_2})^2}{c^{\bfB}}\bigl(1{+} \norm[\L^{\infty}
(\Omega)]{z(\cdot,t)}^2\bigr)\bigr)\norm[\L^{\infty}
(\Omega)]{\theta(\cdot,t)} \bigr)\varphi(t)\\&\textrm{ if }\alpha >0
\end{aligned} 
\end{equation*}
for almost every \((x,t)\in \calQ_{\tau} \).
Thus  $\Theta_{\theta_{\varphi}}$ 
is absolutely continuous on $[0, \tau]$ 
and, also by \eqref{eq:12}, we have
\begin{equation}
\label{eq:15}
\begin{aligned}
&\dot\Theta_{\theta_{\varphi}}(t)= 
\int_{\Omega} c(x)\calG(\theta_{\varphi})\dot\theta_{\varphi}\dd x
= {-}\int_{\Omega}\calG(\theta_{\varphi})
(\dive(\kappa\nabla\theta){+}
\bfA\ee(\dot u){:}\ee(\dot u)\\&
{+}\theta(\beta\bfI{:}\ee(\dot u){+}\partial_zH_2(z){.}\dot z)
{+}
\bfB\dot z{.}\dot z{+}\Psi(\dot z))\varphi \dd x
-\int_{\Omega} c(x) \calG(\theta_{\varphi}) \theta\dot \varphi \dd x
\\&
={-}\int_{\Omega}\calG' (\theta_{\varphi})\kappa\nabla
\theta_{\varphi}{:}\nabla\theta_{\varphi}\dd x 
-\int_{\Omega}\calG(\theta_{\varphi})( 
\bfA\ee(\dot u){:}\ee(\dot u)\\&
{+} \theta (\beta\bfI{:}\ee(\dot u){+}\partial_zH_2(z){.}\dot z) 
{+}\bfB\dot z{.}\dot z {+}\Psi(\dot z)) \varphi \dd x
-\int_{\Omega} c(x) \calG(\theta_{\varphi}) \theta\dot\varphi \dd x
\end{aligned}
\end{equation}
for almost every $t\in[0,\tau]$. We evaluate now the second term
of the right hand side of \eqref{eq:15}. By using
\eqref{eq:H2}, \eqref{eq:LM} and Cauchy-Schwarz's 
inequality, we get
\begin{equation}
\label{eq:13a}
\bfA\ee(\dot u){:}\ee(\dot u)
+\beta\theta\bfI{:}\ee(\dot u)\ge c^{\bfA} 
\abs{\ee(\dot u)}^2 -3\beta\abs{\theta}\abs{\ee(\dot u)} 
\ge \tfrac{c^{\bfA}}{2}\abs{\ee(\dot u)}^2 -\tfrac{ 9 \beta^2}{2 c^{\bfA}}
\abs{\theta}^2,
\end{equation}
and if $\alpha >0$
\begin{equation}
\label{eq:13b}
\begin{aligned}
&\bfB\dot z{.}\dot z + \theta \partial_zH_2(z){.}
\dot z \ge c^{\bfB} \abs{\dot z}^2 - 
\abs{\theta}\abs{\partial_z H_2(z)}\abs{\dot z} 
\\&\geq c^{\bfB} \abs{\dot z}^2 
- C_z^{H_2} \abs{\theta}
(1{+}\abs{z}) \abs{\dot z} 
\geq \tfrac{c^{\bfB}}{2} \abs{\dot z}^2 -
\tfrac{(C_z^{H_2})^2}{c^{\bfB}}(1{+}\abs{z}^{2})\abs{\theta}^2.
\end{aligned}
\end{equation}
We insert \eqref{eq:13a} and \eqref{eq:13b} into \eqref{eq:15}, then
recalling  that $\calG'(\theta_{\varphi}) \ge 0$ 
almost everywhere,
by \eqref{eq:16} and \eqref{eq:16bis}, we obtain
\begin{equation*}
\dot\Theta_{\theta_{\varphi}}(t)\leq\int_{\Omega}\calG(\theta_{\varphi})
 \tfrac{9 \beta^2}{2 c^{\bfA}} \bigl(\abs{\theta}^2 
{+}\tfrac{c(x)}{c^c}\norm[\L^{\infty}(\Omega)]{\theta}\theta\bigr)
\varphi \dd x
\end{equation*}
if $\alpha =0$ and 
\begin{equation*}
\dot\Theta_{\theta_{\varphi}}(t)\leq\int_{\Omega}\calG(\theta_{\varphi})
\bigl( \tfrac{9 \beta^2}{2 c^{\bfA}}{+}
\tfrac{(C_z^{H_2})^2}{c^{\bfB}}  
\bigl(1{+}\norm[\L^{\infty}(\Omega)]{z}^{2 }\bigr)\bigr)\bigl(\abs{\theta}^2 
{+}\tfrac{c(x)}{c^c}\norm[\L^{\infty}(\Omega)]{\theta}\theta\bigr)
\varphi \dd x
\end{equation*}
if $\alpha >0$, for almost every $t\in[0,\tau]$.
Now we observe that \(\calG(\Theta_{\theta_{\varphi}})\)
vanishes whenever $\theta$ is non negative and 
\begin{equation*}
\abs{\theta}^2 
{+}\tfrac{c(x)}{c^c}\norm[\L^{\infty}(\Omega)]{\theta}\theta
= \abs{\theta} \bigl( \abs{\theta}{-} \tfrac{c(x)}{c^c}
\norm[\L^{\infty}(\Omega)]{\theta}\bigr)
\leq \abs{\theta} \bigl( \abs{\theta}{-} 
\norm[\L^{\infty}(\Omega)]{\theta}\bigr) \leq 0
\end{equation*}
whenever $\theta$ is non positive.
Hence \(\dot\Theta_{\theta_{\varphi}}(t)\leq 0\) for almost 
every $t\in[0,\tau]$.
Since \(\Theta_{\theta_{\varphi}}(0)= \int_{\Omega} c(x) 
\Gamma ( - \theta^0(x)) \dd x= 0\), we infer  that
\(\Theta_{\theta_{\varphi}}(t)\leq 0\) for all $t \in [0, \tau]$. 
It follows that $\Gamma (\theta_{\varphi}(x,t) ) =0$ for almost 
every $(x,t) \in \calQ_{\tau}$ implying that 
\(\theta_{\varphi}(x,t) = -\theta(x,t)\varphi(t) \leq 0\) i.e.  
\(\theta(x,t)\geq 0\)
for almost every \((x,t)\in\calQ_{\tau}\). 

Let us  establish now that the temperature \(\theta(x,t)\) 
remains positive for almost every \((x,t)\in\calQ_{\tau}\).
To this aim, we define
\(\tilde \Theta_{\tilde \theta_{\varphi}}(t)
\eqldef\int_{\Omega}c(x)\Gamma(\tilde \theta_{\varphi}(x,t))\dd x\) with
$\tilde \theta_{\varphi}(x,t)\eqldef-\theta(x,t)+\bar\theta\varphi(t)$ 
for almost every $(x,t) \in \calQ_{\tau}$.
Since $\theta\in\W^{1,2}(0, \tau; \L^2(\Omega))$ 
we infer that $\tilde \Theta_{\theta_{\varphi}}$ is
absolutely continuous on $[0,\tau]$ and we have
\begin{equation}
\label{eq:13}
\begin{aligned}
&\dot{\tilde \Theta}_{\tilde \theta_{\varphi}}(t)= 
\int_{\Omega} c(x)\calG(\tilde \theta_{\varphi})\dot{\tilde \theta}_{\varphi}\dd x
= {-}\int_{\Omega}\calG(\tilde \theta_{\varphi})
(\dive(\kappa\nabla\theta){+}
\bfA\ee(\dot u){:}\ee(\dot u)\\&
{+}\theta(\beta\bfI{:}\ee(\dot u){+}\partial_zH_2(z){.}\dot z)
{+}\bfB\dot z{.}\dot z {+}\Psi(\dot z)
{-}c(x)\bar\theta\dot\varphi)\dd x\\&
=-\int_{\Omega}\calG' (\tilde \theta_{\varphi})\kappa\nabla
\tilde \theta_{\varphi}{:}\nabla \tilde \theta_{\varphi}\dd x
{-}\int_{\Omega}\calG(\tilde \theta_{\varphi})( 
\bfA\ee(\dot u){:}\ee(\dot u)\\&
{+} \theta (\beta\bfI{:}\ee(\dot u){+}\partial_zH_2(z){.}\dot z) 
{+}\bfB\dot z{.}\dot z{+}\Psi(\dot z)
{-}c(x)\bar\theta\dot\varphi)\dd x
\end{aligned}
\end{equation}
for almost every $t\in[0,\tau]$.  
We estimate the right hand side of \eqref{eq:13} by using the 
same tricks as previously,
we obtain
\begin{equation*}
\dot{\tilde \Theta}_{\tilde \theta_{\varphi}}(t)\leq
\int_{\Omega}\calG(\tilde \theta_{\varphi})
\bigl( \tfrac{9 \beta^2}{2 c^{\bfA}} \abs{\theta}^2 
{+}c(x)\bar\theta\dot\varphi\bigr) \dd x
\end{equation*}
if $\alpha =0$ and 
\begin{equation*}
\dot{\tilde \Theta}_{\tilde \theta_{\varphi}}(t)
\leq\int_{\Omega}\calG(\tilde \theta_{\varphi})
\bigl(\bigl( \tfrac{9 \beta^2}{2 c^{\bfA}} {+}
\tfrac{(C_z^{H_2})^2}{c^{\bfB}}  
\bigl(1{+}\abs{z}^{2 }\bigr) \bigr)\abs{\theta}^2 
{+}c(x)\bar\theta\dot\varphi\bigr) \dd x
\end{equation*}
if $\alpha >0$, for almost every $t\in[0,\tau]$. It follows 
from \eqref{eq:16} and \eqref{eq:16bis} that
\begin{equation*}
\dot{\tilde \Theta}_{\tilde \theta_{\varphi}} (t)\leq
\int_{\Omega}\calG(\tilde \theta_{\varphi})
\tfrac{9 \beta^2}{2 c^{\bfA}} \bigl(\abs{\theta}^2 
{-}\tfrac{c(x)}{c^c}\norm[\L^{\infty}(\Omega)]{\theta}
\bar\theta\varphi\bigr) \dd x
\end{equation*}
if $\alpha =0$ and
\begin{equation*}
\dot{\tilde \Theta}_{\tilde \theta_{\varphi}} (t)
\leq\int_{\Omega}\calG(\tilde \theta_{\varphi})
\bigl(\tfrac{9 \beta^2}{2 c^{\bfA}}{+}
\tfrac{(C_z^{H_2})^2}{c^{\bfB}}  
\bigl(1{+} \norm[\L^{\infty}(\Omega)]{z}^{2}\bigr)\bigr)\bigl(\abs{\theta}^2 
{-}\tfrac{c(x)}{c^c}\norm[\L^{\infty}(\Omega)]{\theta}
\bar\theta\varphi\bigr) \dd x
\end{equation*}
if $\alpha >0$, for almost every $t\in[0,\tau]$. Then we observe that 
\(\calG(\tilde \theta_{\varphi})\) vanishes whenever 
 \(\theta\geq \bar\theta\varphi\), and 
\begin{equation*}
\abs{\theta}^2 
-\tfrac{c(x)}{c^c}
\norm[\L^{\infty}(\Omega)]{\theta}
\bar\theta\varphi \le \norm[\L^{\infty}(\Omega)]{\theta} 
\bigl(\abs{\theta} {-}\tfrac{c(x)}{c^c} \bar\theta\varphi \bigr) \le
0
\end{equation*}
whenever $0 \le \theta \le \bar \theta \varphi$. 
Since we have already proved that $\theta$ is non 
negative almost everywhere in $\calQ_{\tau}$, we may infer that 
\(\dot{\tilde \Theta}_{\tilde \theta_{\varphi}}(t)\leq 0\)
for almost every $t\in[0,\tau]$.
Therefore 
$\tilde \Theta_{\tilde \theta_{\varphi}}(t)\leq \tilde 
\Theta_{\tilde \theta_{\varphi}}(0)= \int_{\Omega} c(x) 
\Gamma (- \theta^0 {+}\bar \theta) \dd x=0$ 
for all $t\in[0,\tau]$. It follows  
that 
\(\Gamma(\tilde \theta_{\varphi}(x,t) ) = 0\) for almost every  
\((x,t)\in\calQ_{\tau}\),
which implies that
\begin{equation*}
\tilde \theta_{\varphi}(x,t)=-\theta (x,t) +\bar\theta\varphi\leq 0 
\end{equation*}
for almost every \((x,t)\in \calQ_{\tau} \).
\eproof

Furthermore the solutions of problem 
\eqref{eq:ent_eq}--\eqref{eq:init_cond} satisfy the following global estimate:

\begin{proposition}
\label{thm:global-estim}
Assume that \eqref{eq:Psi}, \eqref{eq:H},  \eqref{eq:47}, 
\eqref{eq:LM},  \eqref{eq:Q}, \eqref{eq:f}, \eqref{eq:36},
\(\theta^0\in\W^{1,2}_{\kappa,\rm{Neu}}(\Omega)\), 
\(u^0\in\W^{1,4}_{\rm{Dir}}(\Omega)\)
and \(\zeo\in\W^{2,4}_{\rm{Neu}}(\Omega)\) if \(\alpha>0\) and
\(\zeo\in\L^4(\Omega)\) if \(\alpha=0\) hold. Assume also that 
condition \eqref{init_temp} is satisfied, $\kappa \in 
\C^1({\bar \Omega})$  and $c^{H_1} >0$.
Then, there exists a constant $ C_0>0$,  
depending only 
the data such that
for any solution $(u,z, \theta) $ of problem 
\eqref{eq:ent_eq}--\eqref{eq:init_cond} 
defined on $[0,\tau]$, $\tau \in (0,T]$,  we have
\begin{equation*}
\forall t {\in} [0, \tau]:
\norm[\W^{1,2} (\Omega)]{u (\cdot, t) }^2  + 
\norm[\L^2(\Omega)]{ z (\cdot, t) }^2 + \alpha 
\norm[\L^2(\Omega)]{\nabla z (\cdot, t) }^2
+ \norm[\L^{1}(\Omega)]{\theta (\cdot, t) } \le C_0.
\end{equation*}
\end{proposition}

\bproof
First we choose $\dot u$ as a test-function in \eqref{eq:ent_eq_1} 
and the constant function equal to 1 in \eqref{eq:ent_eq_3}. We get 
\begin{equation}
\label{eqlp:1}
\int_{\calQ_{t}}
(\bfE (\ee (u){-}\bfQ z){+}\beta \theta\bfI{+}\bfA \ee(\dot
u)){:}\ee(\dot u)\dd x\dd s
=\int_{\calQ_{t}} f {\cdot}\dot u\dd x\dd s
\end{equation}
and
\begin{equation}
\label{eqlp:2}
\begin{aligned}
&\int_{\Omega} c (\cdot) \theta(\cdot,t)\dd x=
\int_{\Omega} c (\cdot) \theta^0\dd x
+\int_{\calQ_{t}}\bfA \ee(\dot u){:}\ee(\dot u)\dd x\dd s
+\int_{\calQ_{t}}\bfB \dot z{.}\dot z\dd x\dd s\\&
+\int_{\calQ_{t}}\theta(\beta \bfI {:} \ee(\dot u){+} \partial_z
H_2(z){.}
\dot z)\dd x\dd s+
\int_{\calQ_{t}}\Psi(\dot z)\dd x\dd s.
\end{aligned}
\end{equation}
Then we use the definition of the subdifferential $\partial \Psi$; 
for almost every $s \in [0, \tau]$ and all 
$\tilde z \in \L^2(\Omega, {\mathcal Z})$, we have  
\begin{equation*}
\begin{aligned}
& \int_{\Omega}
(\bfB\dot z (\cdot, s) {-} {\tilde \bfQ}^{\tra}
\bfE (\ee(u(\cdot, s) ){-} \bfQ z(\cdot, s) ){.}
(\tilde z {-} \dot  z(\cdot, s) )\dd x\\&
{-}  \int_{\Omega} \alpha \Delta z(\cdot, s) {.}
(\tilde z {-} \dot  z(\cdot, s) ) \dd x \\
&
{+} \int_{\Omega}  \partial_z H_1(z(\cdot, s) ){+}
\theta (\cdot, s) \partial_zH_2(z(\cdot, s) ) {.}
(\tilde z {-} \dot  z(\cdot, s) ) \dd x  \\
&
+\int_{\Omega}\Psi(\tilde z) \dd x -\int_{\Omega}
\Psi(\dot  z(\cdot, s) ) \dd x \ge 0.
\end{aligned}
\end{equation*}
But $\Psi$ is positively homogeneous of degree 1, so by 
choosing successively $\tilde z \equiv 0$ and 
$\tilde z= 2 \dot z (\cdot, s)$ and integrating over 
$[0,t] \subset [0, \tau]$, we obtain
\begin{equation} \label{eqlp:3}
\begin{aligned}
& \int_{\calQ_t}
(\bfB\dot z  {-} {\tilde \bfQ}^{\tra}\bfE (\ee(u ){-} 
\bfQ z ){+}\partial_z H_1(z ){+}
\theta (\cdot, s) \partial_zH_2(z ){-}\alpha \Delta z ){.}
 \dot  z  \dd x \dd s \\
 &
+\int_{\calQ_t}\Psi(\dot  z) \dd x \dd s  =  0.
\end{aligned}
\end{equation}
Now we add \eqref{eqlp:1}, \eqref{eqlp:2} and \eqref{eqlp:3}, we get 
\begin{equation}
\label{eqlp:4}
\begin{aligned}
&\tfrac12\int_{\Omega}\bfE (\ee(u(\cdot,t)){-}\bfQ z(\cdot,t)){:}
(\ee(u(\cdot,t)){-} \bfQ z(\cdot,t))\dd x
+ \int_{\Omega} H_1(z(\cdot,t))\dd x 
\\&
+\tfrac{\alpha}2\norm[\L^2(\Omega)]{\nabla z(\cdot, t)}^2
 + \int_{\Omega} c (\cdot) \theta(\cdot,t)\dd x\\&
= \tfrac12\int_{\Omega}\bfE (\ee(u^0){-}\bfQ z^0){:}
(\ee(u^0){-}\bfQ z^0 )\dd x \\
&
+
\tfrac{\alpha}2\norm[\L^2(\Omega)]{\nabla z^0}^2 
 + \int_{\Omega} 
H_1(z^0) \dd x
+ \int_{\Omega} c (\cdot) \theta^0\dd x
+\int_{\calQ_{t}}f {\cdot}\dot u\dd x\dd s.
\end{aligned}
\end{equation}
We estimate from below the two   first   terms of 
the left hand side by using \eqref{eq:H1}, \eqref{eq:3} 
and \eqref{eq:Q}. Indeed, for any $\lambda \in (0,1)$, we find
\begin{equation*}
\begin{aligned}
& \tfrac{1}{2} \int_{\Omega}\bfE (\ee(u(\cdot,t)){-}\bfQ z(\cdot,t))
{:}(\ee(u(\cdot,t)){-} \bfQ z(\cdot,t))\dd x +   
\int_{\Omega} H_1(z(\cdot,t))\dd x 
 \\ &
\ge \tfrac{1}{2} (1{-} \lambda) c^{\bfE} 
\norm[\L^2(\Omega)]{\ee (u (\cdot, t) )}^2
  +  \bigl( 1{-} \tfrac{1}{\lambda} \bigr)  
\norm[\L^{\infty}(\Omega)]{\bfE} 
\bigl( \norm{\tilde \bfQ}^2 \norm[\L^2(\Omega)]{z (\cdot, t) }^2 
{+} \abs{\mathrm{Q}}^2 \abs{\Omega} \bigr) \\
&  +  c^{H_1} \norm[\L^2(\Omega)]{z (\cdot, t) }^2 - \tilde c^{H_1} \abs{\Omega}.
\end{aligned}
\end{equation*}
 We may choose $\lambda \in (0,1)$ such that 
\begin{equation*}
 \bigl( 1{-} \tfrac{1}{\lambda} \bigr)  
\norm[\L^{\infty}(\Omega)]{\bfE}  \norm{\tilde \bfQ}^2  + c^{H_1} >0,
\end{equation*}
i.e.
\begin{equation*}
1> \lambda > \tfrac{ \norm[\L^{\infty}(\Omega)]{\bfE}  
\norm{\tilde \bfQ}^2}{ \norm[\L^{\infty}(\Omega)]{\bfE}  
\norm{\tilde \bfQ}^2 + c^{H_1}}.
\end{equation*}
Then
\begin{equation*}
\begin{aligned}
& \tfrac{1}{2} \int_{\Omega}\bfE (\ee(u(\cdot,t)){-}\bfQ z(\cdot,t)){:}
(\ee(u(\cdot,t)){-} \bfQ z(\cdot,t))\dd x +   \int_{\Omega} H_1(z(\cdot,t))\dd x 
 \\ &
\ge C \bigl( \norm[\W^{1,2} (\Omega)]{u(\cdot,t)}^2 {+} 
\norm[\L^2(\Omega)]{ z(\cdot,t)}^2 \bigr) - \tilde C,
\end{aligned}
\end{equation*}
with
\begin{equation*}
\begin{aligned}
& C\eqldef\min \bigl( \tfrac{1}{2} (1{-} \lambda) c^{\bfE} \cK ,  
\bigl( 1{-} \tfrac{1}{\lambda} \bigr)  
\norm[\L^{\infty}(\Omega)]{\bfE}  \norm{\tilde \bfQ}^2  {+} c^{H_1} \bigr) , \\
& \tilde C \eqldef \bigl(\tfrac{1}{\lambda} {-}1 \bigr)  
\norm[\L^{\infty}(\Omega)]{\bfE} \abs{\mathrm{Q}}^2 \abs{\Omega} 
+ \tilde c^{H_1} \abs{\Omega}.
\end{aligned}
\end{equation*}
Now we integrate by parts the last term of the right hand side of 
\eqref{eqlp:4} to get
\begin{equation*}
\begin{aligned}
& \tfrac{C}{2}  \norm[\W^{1,2} (\Omega)]{u(\cdot,t)}^2 + C 
\norm[\L^2(\Omega)]{ z(\cdot,t)}^2 
+\tfrac{\alpha}2\norm[\L^2(\Omega)]{\nabla z(\cdot,t)}^2  +
\int_{\Omega} c (\cdot) \theta (\cdot, t) \dd x \\
& \leq  \tfrac12\int_{\Omega}\bfE (\ee(u^0){-}\bfQ z^0){:}
(\ee(u^0){-}\bfQ z^0)\dd x+
\tfrac{\alpha}2\norm[\L^2(\Omega)]{\nabla z^0}^2 \\
& + \int_{\Omega} 
H_1(z^0) \dd x
+ \int_{\Omega} c(\cdot) \theta^0\dd x + \tilde C
+\norm[{\C^0([0,T]; \L^2(\Omega))}]{f} \norm[\L^2(\Omega)]{u^0} \\
& +\tfrac{1}{2 C} \norm[{\C^0([0,T]; \L^2(\Omega))}]{f}^2 
+ \tfrac12 \norm[\L^2(0,T; \L^2(\Omega))]{\dot f}^2  
+\tfrac12 \int_0^{t} \norm[\L^2(\Omega)]{u}^2 \dd s.
\end{aligned}
\end{equation*}
Then, recalling that $\theta$ remains non negative, Gr\"onwall's 
lemma allows us to conclude.
\eproof

Let us find now some sufficient conditions on the data which will lead to a 
global existence result, i.e. 
existence of a solution of problem \eqref{eq:ent_eq}--\eqref{eq:init_cond}  
defined on the whole interval $[0,T]$.
First we observe that the heat-transfer equation \eqref{eq:ent_eq_3} and 
the system composed of the momentum equilibrium equation and the flow 
rule \eqref{eq:ent_eq_1}--\eqref{eq:ent_eq_2} are totally decoupled if 
$\beta =0$ and $\partial_{\ze} H_2 \equiv 0$. In such a case, we may obtain 
a solution of \eqref{eq:ent_eq}--\eqref{eq:init_cond} by applying 
Proposition \ref{sec:existence} to solve 
\eqref{eq:ent_eq_1_1}--\eqref{eq:init_cond2} with $\tilde \theta =0$ and 
$\tau =T$, then by finding the solution  $\theta$ of 
\eqref{eq:19}--\eqref{eq:21} with
\begin{equation*}
f^{\tilde \theta} = \bfA\ee(\dot u){:}\ee(\dot u)
+ \bfB\dze{.}\dze +\Psi(\dze).
\end{equation*}
Hence we will consider only the case of non vanishing coupling 
parameters $\beta \not=0$ or $\partial_{\ze} H_2 \not\equiv 0$. By 
using more detailed estimates for the mapping $\tilde \theta 
\mapsto (u,z)$, we can obtain more precise estimates for $f^{\tilde \theta}$ 
which will allow us to prove that the mapping $\Phi^{\tilde \theta, 
\theta}_{T} $ possesses a fixed point in $\C^0([0,T]; \L^4(\Omega))$. 

Let us begin with the case $\alpha =0$. Then we have 

\begin{lemma}({\cite[Thm.~4.1]{PaoPet11GETV}}). 
\label{lemma1}
Let $\tau \in (0,T]$. Assume that \eqref{eq:Psi}, 
\eqref{eq:H},  \eqref{eq:47}, \eqref{eq:LM}, \eqref{eq:Q}, \eqref{eq:f} hold.
Let $\tilde  \theta \in \L^q(0,\tau; \L^p(\Omega))$, with $q >8 $ and  
$p \in [4,6]$,  $u^0 \in \W^{1,p}_{\rm{Dir}}(\Omega)$ and 
$z^0 \in \L^p(\Omega)$ be given and denote by $(u,z)$ the 
unique solution of \eqref{eq:ent_eq_1_1}--\eqref{eq:init_cond2}. 
Then, there exists a  non decreasing positive mapping 
$\tau \mapsto C_{u,z}(\tau)$, independent of the initial data,  such that
\begin{equation*}
\begin{aligned}
&\|u\|_{\C^0([0,\tau]; \W^{1,p}(\Omega))} + \|z\|_{\C^0([0,\tau]; 
\L^p(\Omega))} + \norm[\L^q(0,\tau; \W^{1,p}(\Omega))]{\dot  u} 
+\norm[\L^q(0,\tau; \L^p(\Omega))]{\dot   z} 
\\&\leq C^q_{u,z}(\tau)  
\bigl( \norm[\W^{1,p}(\Omega)]{u^0} {+} \norm[\L^p(\Omega)]{z^0} 
{+}\beta \norm[\L^q(0,\tau; \L^p(\Omega))]{\tilde \theta}{+}1 \bigr).
\end{aligned}
\end{equation*}
\end{lemma}

Let us assume from now on that $u^0 \in \W^{1,4}_{\textrm{Dir}}(\Omega)$, 
$z^0 \in \L^4(\Omega)$, $\theta^0 \in \W^{1,2}_{\textrm{$\kappa$,Neu}}(\Omega)$ 
and let $\tilde \theta \in \C^0 ([0, \tau]; \L^4(\Omega))$ with 
$\tau \in (0, T]$.  From Lemma \ref{lemma1} we can estimate 
$f^{\tilde \theta}$ as follows
\begin{equation*}
\begin{aligned}
& \norm[\L^{q/2}(0, \tau; \L^2(\Omega))]{f^{\tilde \theta}} \\&
\le \norm{\bfA} \norm[\L^q (0, \tau; \L^4(\Omega))]{\ee( \dot u)}^2 
+ 3 \beta \norm[\L^q (0, \tau; \L^4(\Omega))]{\tilde \theta} 
\norm[\L^q (0, \tau; \L^4(\Omega))]{\ee( \dot u)} \\
&  + \norm{\bfB} \norm[\L^q (0, \tau; \L^4(\Omega))]{\dot \ze}^2 
+ C^{\Psi} \norm[\L^{q/2} (0, \tau; \L^2(\Omega))]{\dot \ze} \\
& \le \bigl( \norm{\bfA} {+} \norm{\bfB}\bigr) (C^q_{u,z}(\tau))^2 
\bigl( \norm[\W^{1,4}(\Omega)]{u^0} {+} \norm[\L^4(\Omega)]{z^0} 
{+}\beta \norm[\L^q(0,\tau; \L^4(\Omega))]{\tilde \theta}{+}1 \bigr)^2 
\\& {+} \bigl( 3 \beta   \norm[\L^q(0,\tau; \L^4(\Omega))]{\tilde \theta} 
{+} C^{\Psi} \abs{\Omega}^{\frac{1}{4}} \tau^{\frac{1}{q}}\bigr)\\& 
C^q_{u,z}(\tau) \bigl( \norm[\W^{1,4}(\Omega)]{u^0} {+} \norm[\L^4(\Omega)]{z^0} 
{+}\beta \norm[\L^q(0,\tau; \L^4(\Omega))]{\tilde \theta}{+}1 \bigr)
\\ &
\le C^q_{f^{\tilde \theta}} (\tau,  \norm[\W^{1,4}(\Omega)]{u^0}, 
\norm[\L^4(\Omega)]{z^0} ) 
\bigl( 1 {+} \beta^2 \norm[\L^q(0,\tau; \L^4(\Omega))]{\tilde \theta}^2 \bigr), 
\end{aligned}
\end{equation*}
where $C^q_{f^{\tilde \theta}}\bigl(\tau,  \norm[\W^{1,4}(\Omega)]{u^0}, 
\norm[\L^4(\Omega)]{z^0}\bigr) $ is given by
\begin{equation*}
\begin{aligned}
& C^q_{f^{\tilde \theta}}  \bigl(\tau,  \norm[\W^{1,4}(\Omega)]{u^0}, 
\norm[\L^4(\Omega)]{z^0} \bigr)\eqldef
\max \bigl(  
2 \bigl(\norm{\bfA} {+}\norm{\bfB}\bigr) 
(C^q_{u,z} (\tau))^2{+} 3 C^q_{u,z} (\tau) {+} 1 , \\
&  2 \bigl( \norm{\bfA} {+} \norm{\bfB} {+}  \tfrac{9}{4} \bigr) 
(C^q_{u,z} (\tau))^2 
\bigl( \norm[\W^{1,4}(\Omega)]{u^0} {+} \norm[\L^4(\Omega)]{z^0} {+}1 \bigr)^2 
\\&
{+} \tfrac{1}{2} \bigl(C^{\Psi} \abs{\Omega}^{\frac{1}{4}} \tau^{\frac{1}{q}}
C^q_{u,z}(\tau)\bigr)^2
{+}  C^{\Psi} \abs{\Omega}^{\frac{1}{4}} \tau^{\frac{1}{q}} 
C^q_{u,z} (\tau) \bigl( \norm[\W^{1,4}(\Omega)]{u^0} {+} 
\norm[\L^4(\Omega)]{z^0} {+}1 \bigr) 
\bigr)
\end{aligned}
\end{equation*}
for any $q>8$. It follows  that 
$\theta = \Phi^{\tilde \theta, \theta}_{\tau} (\tilde \theta)$ 
can be estimated as
\begin{equation*} 
\begin{aligned}
\norm[\L^{\infty}(0, \tau; \W^{1,2} (\Omega))]{\theta} 
&\leq C_{\theta} \exp(\tfrac{\tau}{c^c})  \bigl(  
\norm[\W^{1,2}(\Omega)]{\theta^0} {+}  \norm[\L^{2} (0,\tau; 
\L^2(\Omega))]{f^{\tilde \theta}} \bigr) \\
& \le C_{\theta} \exp(\tfrac{\tau}{c^c})  \bigl(  
\norm[\W^{1,2}(\Omega)]{\theta^0} {+}  \tau^{\frac{q-4}{2q}}
\norm[\L^{q/2} (0,\tau; \L^2(\Omega))]{f^{\tilde \theta}} \bigr)
\end{aligned}
\end{equation*}
where $C_{\theta}$ is the constant, independent of $\tau$ and of 
the initial data, introduced in Proposition \ref{prop:Phi1} (see \eqref{eq:9}). 
Since $\theta \in \C^0([0, \tau]; \L^4(\Omega))$ and 
$\W^{1,2} (\Omega) \hookrightarrow \L^4(\Omega)$, we obtain
\begin{equation*}
\begin{aligned}
& \|\theta\|_{C^0([0, \tau]; \L^4(\Omega))} \le  C_1 C_{\theta}   
\exp(\tfrac{\tau}{c^c}) \bigl( \norm[\W^{1,2} (\Omega)]{ \theta^0}  \\
 & +  \tau^{\frac{q-4}{2q}} C^q_{f^{\tilde \theta}}  \bigl(\tau,  
\norm[\W^{1,4}(\Omega)]{u^0}, \norm[\L^4(\Omega)]{z^0} \bigr) 
\bigl( 1 {+} \beta^2 
\norm[\L^q(0,\tau; \L^4(\Omega)]{\tilde \theta}^2 \bigr) \\
& \le 
 C^q
\bigl(\tau, \norm[\W^{1,2}(\Omega)]{\theta^0}, \norm[\W^{1,4}(\Omega)]{u^0}, 
\norm[\L^4(\Omega)]{z^0}\bigr) \bigl( 1 {+} \beta^2 \tau^{\frac{2}{q}} \| 
\tilde \theta\|_{\C^0([0, \tau]; \L^4(\Omega))}^2 \bigr), 
\end{aligned}
\end{equation*}
where $C_1$ is the generic constant involved in the continuous 
embedding of  $\W^{1,2} (\Omega) $ into $ \L^4(\Omega)$ and 
\begin{equation*}
\begin{aligned}
& C^q
\bigl(\tau, \norm[\W^{1,2}(\Omega)]{\theta^0}, \norm[\W^{1,4}(\Omega)]{u^0}, 
\norm[\L^4(\Omega)]{z^0}\bigr) \\
& 
\eqldef  
C_1 C_{\theta}   \exp(\tfrac{\tau}{c^c}) 
\bigl(\norm[\W^{1,2} (\Omega)]{\theta^0} +  \tau^{\frac{q-4}{2q}} 
C^q_{f^{\tilde \theta}}  \bigl( \tau ,  \norm[\W^{1,4}(\Omega)]{u^0}, 
\norm[\L^4(\Omega)]{z^0} \bigr)\bigr).
 \end{aligned}
\end{equation*}

We cannot expect to get a global existence result without further 
assumptions on $\beta$. This is not very surprising since 
$f^{\tilde \theta}$ behaves as a quadratic coupling term if 
$\beta>0$. But the mapping 
\begin{equation*}
\gamma^q : R^{\theta} 
\mapsto C^q\bigl(T, \norm[\W^{1,2}(\Omega)]{\theta^0}, 
\norm[\W^{1,4}(\Omega)]{u^0}, \norm[\L^4(\Omega)]{z^0}\bigr)  
\bigl(1 {+} \beta^2  T^{\frac{2}{q}} (R^{\theta})^2\bigr) - R^{\theta}
\end{equation*}
admits a minimum for $R^{\theta}= R^{\theta}_{q,\textrm{min}} \eqldef 
\tfrac{1}{2 C^q(T, \norm[\W^{1,2}(\Omega)]{\theta^0}, 
\norm[\W^{1,4}(\Omega)]{u^0}, \norm[\L^4(\Omega)]{z^0})  
\beta^2  T^{\frac{2}{q}}}$ and 
\begin{equation*}
\gamma^q (R^{\theta}_{q,\textrm{min}}) = C^q\bigl(T, 
\norm[\W^{1,2}(\Omega)]{\theta^0}, \norm[\W^{1,4}(\Omega)]{u^0}, 
\norm[\L^4(\Omega)]{z^0}\bigr)  - \tfrac{R^{\theta}_{q,\textrm{min}}}{2} .
\end{equation*}
Hence $\gamma^q (R^{\theta}_{q,\textrm{min}}) <0$ if 
$R^{\theta}_{q, \textrm{min}} >2  
C^q\bigl(T, \norm[\W^{1,2}(\Omega)]{\theta^0}, \norm[\W^{1,4}(\Omega)]{u^0}, 
\norm[\L^4(\Omega)]{z^0}\bigr)   $, i.e. 
\begin{equation} 
\label{eq:cs1}
0 <  \beta < \tfrac{1}{2  
C^q(T, \norm[\W^{1,2}(\Omega)]{\theta^0}, 
\norm[\W^{1,4}(\Omega)]{u^0}, \norm[\L^4(\Omega)]{z^0})  T^{\frac{1}{q}}}.
\end{equation}

Let us fix now $q>8$ and assume  that this condition on $\beta$ holds. 
We  choose $R^{\theta}= R^{\theta}_{q,\textrm{min}}$. We may observe that, 
since $\beta$ satisfies condition \eqref{eq:cs1}, we have
\begin{equation*}
\begin{aligned}
 R^{\theta}_{q, \textrm{min}} & = 
\tfrac{1}{2 C^q(T, \norm[\W^{1,2}(\Omega)]{\theta^0}, 
\norm[\W^{1,4}(\Omega)]{u^0}, \norm[\L^4(\Omega)]{z^0})  \beta^2  
T^{\frac{2}{q}}} \\
& > 2 C^q\bigl(T, \norm[\W^{1,2}(\Omega)]{\theta^0}, 
\norm[\W^{1,4}(\Omega)]{u^0}, 
\norm[\L^4(\Omega)]{z^0}\bigr)  
\\&> C_1 C_{\theta}   \exp\bigl(\tfrac{T}{c^c}\bigr)  
\norm[\W^{1,2} (\Omega)]{\theta^0}.
\end{aligned}
\end{equation*}
Thus we can apply the results of Corollary \ref{corollaire}: there 
exists $\tau \in (0,T]$ such that $\Phi^{\tilde \theta, \theta}_{\tau} $ 
possesses a fixed point in $\C^0([0,\tau]; \L^4(\Omega))$.  But the 
previous estimate implies also that 
\begin{equation*}
\begin{aligned}
& \|\Phi^{\tilde \theta, \theta}_{\tau} (\tilde \theta) \|_{\C^0([0, \tau]; 
\L^4(\Omega))} = \| \theta  \|_{\C^0([0, \tau]; \L^4(\Omega))} \\
& \le C^q\bigl(\tau, \norm[\W^{1,2}(\Omega)]{\theta^0}, 
\norm[\W^{1,4}(\Omega)]{u^0}, \norm[\L^4(\Omega)]{z^0}\bigr)  
\bigl( 1 {+} \beta^2 \tau^{\frac{2}{q}} \| \tilde  \theta  \|^2_{\C^0([0,
  \tau]; 
\L^4(\Omega))}\bigr) \\
& \le C^q\bigl(T, \norm[\W^{1,2}(\Omega)]{\theta^0}, 
\norm[\W^{1,4}(\Omega)]{u^0}, 
\norm[\L^4(\Omega)]{z^0}\bigr)\bigl(1{+} 
\beta^2 T^{\frac{2}{q}}( R^{\theta}_{q, \textrm{min}})^2\bigr) \\
& = \gamma^q(R^{\theta}_{q, \textrm{min}}) + R^{\theta}_{q, \textrm{min}}   
< R^{\theta}_{q, \textrm{min}}
\end{aligned}
\end{equation*}
for any $\tau \in (0,T]$ and any $\tilde \theta 
\in {\bar B}_{\C^0([0, \tau]; \L^4(\Omega))}(0, R^{\theta}_{q,\textrm{min}})$.
Hence we can consider $\tau =T$ and the closed convex bounded set 
${\mathcal C} \eqldef {\bar B}_{\C^0([0, T]; \L^4(\Omega))} 
(0, R^{\theta}_{q, \textrm{min}})$. We have $\Phi^{\tilde \theta, \theta}_{T} 
({\mathcal C}) \subset {\mathcal C}$, and using Schauder's 
fixed point theorem, we infer that $\Phi^{\tilde \theta, \theta}_{T}$ 
admits a fixed point $\theta$ in $\C^0([0, T]; \L^4(\Omega))$. 
Then we define $(u, \ze)$ as the unique solution of 
\eqref{eq:ent_eq_1_1}--\eqref{eq:init_cond2} with $\tilde 
\theta = \theta$ and $\tau=T$. By definition of 
$\Phi^{\tilde \theta, \theta}_{T} $, 
$(u,z, \theta)$ is a global solution of the coupled problem 
\eqref{eq:ent_eq}--\eqref{eq:init_cond} on $[0,T]$.

Now let us consider the case $\alpha >0$. 

\begin{lemma}
({\cite[Lemma~4.4]{PaoPet11GESM} and 
\cite[Lemma~3.4]{PaoPet11TMPV}}). 
\label{lemma2}
Let $\tau \in (0,T]$. Assume that \eqref{eq:Psi}, \eqref{eq:H},  
\eqref{eq:47}, \eqref{eq:LM}, \eqref{eq:Q}, \eqref{eq:f} hold.
Let $\tilde  \theta \in \L^q(0,\tau; \L^p(\Omega))$, with $q >8 $ 
and  $p \in [4,6]$,  $u^0 \in \W^{1,p}_{\rm{Dir}}(\Omega)$ and 
$z^0 \in \W^{2,p}_{\rm{Neu}} (\Omega)$ be given and denote by 
$(u,z)$ the unique solution of \eqref{eq:ent_eq_1_1}--\eqref{eq:init_cond2}.
Then, there exists a non-decreasing positive mapping 
$\tau \mapsto C^q_{u}(\tau)$, independent of the initial data,  such that
\begin{equation*}
\begin{aligned}
&\|u\|_{\C^0([0, \tau]; \W^{1,p}(\Omega))} + 
\norm[\L^q(0,\tau;\W^{1,p} (\Omega))]{\dot u}
\\&\leq
C^q_{u} (\tau) \bigl(
\norm[\L^{q} (0,\tau;\W^{1,2} (\Omega))]{z}{+}\beta 
\norm[\L^{q} (0,\tau;\L^p(\Omega))]{\tilde\theta} 
{+}\norm[\W^{1,p}(\Omega)]{u^0}{+}1\bigr).
\end{aligned}
\end{equation*}
\end{lemma}

Let $u^0 \in \W^{1,4}_{\textrm{Dir}}(\Omega)$,
$z^0 \in \W^{2,4}_{\textrm{Neu}} (\Omega)$ and $\theta^0 
\in \W^{1,2}_{\textrm{$\kappa$, Neu}} (\Omega)$ and let 
$\tilde \theta \in \C^0 ([0, \tau]; \L^4(\Omega))$ 
with $\tau \in (0, T]$. With similar computations as 
in Lemma \ref{lm:cont_uz} and Proposition \ref{thm:global-estim}, 
we can obtain

\begin{lemma}  
\label{lemma3}
Let $\tau \in (0, T]$. Assume that \eqref{eq:Psi}, \eqref{eq:H},  
\eqref{eq:47}, \eqref{eq:LM}, \eqref{eq:Q}, \eqref{eq:f}, 
$u^0 \in \W^{1,4}_{\rm{Dir}}(\Omega)$ and $z^0 \in \W^{2,4}_{\rm{Neu}} 
(\Omega)$ hold. Let $\tilde \theta \in \C^0 ([0, \tau]; \L^4(\Omega))$ 
be given and denote by $(u,z)$  the unique solution of 
\eqref{eq:ent_eq_1_1}--\eqref{eq:init_cond2}. Then 
\begin{equation*}
\begin{aligned} 
&\norm[\L^{\infty}(0, \tau; \W^{1,2} (\Omega))]{u}^2  + 
\norm[\L^{\infty}(0, \tau; \W^{1,2}(\Omega))]{ z}^2 
\\&
\leq C \bigl(\norm[\W^{1,2}(\Omega)]{u^0},  \norm[\W^{1,2}(\Omega)]{z^0} \bigr)  
(X{+}1) \exp(c_0(X{+}1)\tau),
\end{aligned}
\end{equation*}
where $X\eqldef \bigl(\beta^2 {+} (C_{\ze}^{H_2})^2\bigr) \| \tilde \theta 
\|^2_{\C^0 ([0, \tau]; \L^4(\Omega))}$, $c_0>0$ is a constant 
independent of the initial data and $\tau$, and 
$C \bigl(\norm[\W^{1,2}(\Omega)]{u^0},  \norm[\W^{1,2}(\Omega)]{z^0}\bigr)$ 
is a non decreasing positive function  of each of its arguments.
\end{lemma}

\bproof
Let $C^{H_1}>0$ and  define
\begin{equation*} 
\delta(t)\eqldef
\tfrac12 \int_{\Omega}
\bfE(\ee(  u ){-} \bfQ  \ze ){:}
(\ee(  u ){-}  \bfQ \ze )
\dd x
- 
\tfrac{\alpha}2 \int_{\Omega}\Delta \ze{.}\ze
\dd x+ \tfrac{C^{H_1}}{2} \int_{\Omega} \abs{\ze}^2 \dd x
\end{equation*}
for all \(t \in [0,\tau]\). 
As in Proposition \ref{thm:global-estim} we can check that, 
for any $\lambda \in (0,1)$, we have
\begin{equation*}
\begin{aligned}
& \tfrac{1}{2} \int_{\Omega}\bfE (\ee(u){-}\bfQ z){:}(\ee(u )
{-} \bfQ z )\dd x +   \tfrac{C^{H_1}}{2} \int_{\Omega} \abs{\ze}^2 \dd x
\ge \tfrac{1}{2} (1{-}\lambda) c^{\bfE} \norm[\L^2(\Omega)]{\ee (u)}^2\\&
+ \bigl(1{-} \tfrac{1}{\lambda} \bigr)  \norm[\L^{\infty}(\Omega)]{\bfE} 
\bigl( \norm{\tilde \bfQ}^2 \norm[\L^2(\Omega)]{z}^2 {+} 
\abs{{\mathrm{Q}}}^2 
\abs{\Omega} \bigr)
+\tfrac{C^{H_1}}{2} \norm[\L^2(\Omega)]{z}^2.
\end{aligned}
\end{equation*}
Thus we may choose $\lambda \in (0,1)$ such that 
\begin{equation*}
1> \lambda > \tfrac{ \norm[\L^{\infty}(\Omega)]{\bfE}  
\norm{\tilde \bfQ}^2}{ \norm[\L^{\infty}(\Omega)]{\bfE}  
\norm{\tilde \bfQ}^2 + \frac{C^{H_1}}{2} },
\end{equation*}
and we obtain 
\begin{equation*}
\delta (t)
\ge C_{\delta} \bigl( \norm[\W^{1,2} (\Omega)]{u(\cdot,t)}^2 
{+} \norm[\W^{1,2} (\Omega)]{ z(\cdot,t)}^2  \bigr)  - \tilde C_{\delta}
\end{equation*}
for all \(t \in [0,\tau]\), with
\begin{equation*}
\begin{aligned}
&C_{\delta} \eqldef\min \bigl( \tfrac{1}{2} (1{-} \lambda) c^{\bfE} \cK ,  
\bigl( 1{-} \tfrac{1}{\lambda} \bigr)  \norm[\L^{\infty}(\Omega)]{\bfE}  
\norm{\tilde \bfQ}^2 {+} \tfrac{C^{H_1} }{2}, \tfrac{\alpha}{2} \bigr)\\& 
\tilde C_{\delta}\eqldef  \bigl(\tfrac{1}{\lambda}{-}1 \bigr)  
\abs{\mathrm{Q}}^2 
\abs{\Omega}.
\end{aligned}
\end{equation*}
Moreover $\delta$ is absolutely continuous on $[0, \tau]$ and, 
by similar computations as in Lemma \ref{lm:cont_uz}, we get
\begin{equation*}
\begin{aligned}
&
\dot\delta (t)
{+}c^{\bfA}
\norm[\L^2(\Omega)]{\ee(\dot{ u})}^2 
{+}c^{\bfB}
\norm[\L^2(\Omega)]{\dot{\ze}}^2
\le 
  C^{H_1}  \int_{\Omega}  \ze {.} \dot\ze\dd x 
{-}\beta\int_{\Omega}\tilde \theta
\bfI{:}\ee(\dot{ u})\dd x \\
&
{-}\int_{\Omega}  \partial_{\ze}H_1(\ze){.}\dot{ \ze}\dd x
{-}\int_{\Omega} \tilde \theta \partial_{\ze}H_2(\ze){.}
\dot{ \ze}\dd x {+} \int_{\Omega} f {.} \dot u \dd x 
\end{aligned}
\end{equation*}
for almost every \(t\in [0,\tau]\).
We estimate the right hand side of this last inequality by using
\eqref{eq:H2},
we obtain 
\begin{equation*}
\begin{aligned}
& \dot\delta (t)
+c^{\bfA}
\norm[\L^2(\Omega)]{\ee(\dot{ u})}^2 
+c^{\bfB}
\norm[\L^2(\Omega)]{\dot{\ze}}^2 
\leq
(C^{H_1}{+}C_{\ze}^{H_1})  \int_{\Omega} \abs{ \ze} \abs{\dot{ \ze}} \dd x
+ C_{\ze}^{H_1} \int_{\Omega}  \abs{\dot{ \ze}} \dd x \\
& + 3 \beta \int_{\Omega} \abs{\tilde \theta}
\abs{\ee(\dot{ u})} \dd x + C_{\ze}^{H_2}    \int_{\Omega} 
\abs{\tilde \theta} (1{+} \abs{ \ze}) \abs{\dot{ \ze}} \dd x 
+ \int_{\Omega} \abs{f} \abs{ \dot u} \dd x.
\end{aligned}
\end{equation*}
Then, with  Cauchy-Schwarz's inequality
\begin{equation*}
\begin{aligned}
&\dot \delta (t) 
+\tfrac{c^{\bfA}}2
\norm[\L^2(\Omega)]{\ee(\dot{ u})}^2 
+\tfrac{c^{\bfB}}4
\norm[\L^2(\Omega)]{\dot{ \ze}}^2
\\&\leq \bigl( \tfrac{(C^{H_1}{+} C_{\ze}^{H_1}  )^2}{c^{\bfB}} 
{+} \tfrac{2C_1^2 }{c^{\bfB}} (C_z^{H_2})^2 \| \tilde \theta \|^2_{\C^0([0,
  \tau]; 
\L^4(\Omega))} \bigr) \norm[\W^{1,2}(\Omega)]{z (\cdot, t) }^2 \\
& + \tfrac{ (C_z^{H_1})^2}{c^{\bfB}} \abs{\Omega} 
+ \bigl( \tfrac{9\beta^2}{c^{\bfA}} {+} \tfrac{2 (C_z^{H_2})^2}{c^{\bfB}} \bigr) 
\abs{\Omega}^{\frac{1}{2}}  \| \tilde \theta \|^2_{\C^0([0, \tau]; \L^4(\Omega))} 
+ \tfrac{1}{c^{\bfA} \cK} \norm[\L^{\infty} (0,T; \L^2(\Omega))]{f}^2
\end{aligned}
\end{equation*}
for almost every \(t\in [0,\tau]\), where we recall that 
$C_1$ is the generic constant involved in the continuous 
embedding of $\W^{1,2}(\Omega)$ into $\L^4(\Omega)$. Since 
\begin{equation*}
\norm[\W^{1,2}(\Omega)]{z (\cdot, t) }^2 \le \tfrac{ \delta (t) {+} 
\tilde C_{\delta}}{C_{\delta}}
\end{equation*}
for all $t \in [0,\tau]$, we may define $c_0$ and 
$C\bigl(\norm[\W^{1,2}(\Omega)]{u^0},  \norm[\W^{1,2}(\Omega)]{z^0} \bigr)$ by
\begin{equation*}
\begin{aligned}
 c_0\eqldef& \tfrac{1}{C_{\delta}} \max \bigl( \tfrac{
(C^{H_1}{+} C_{\ze}^{H_1}  )^2}{c^{\bfB}}, \tfrac{2C_1^2 }{c^{\bfB}}  \bigr), \\
C\bigl(\norm[\W^{1,2}(\Omega)]{u^0},  \norm[\W^{1,2}(\Omega)]{z^0}\bigr) 
\eqldef& \tfrac{\delta (0) + \tilde C_{\delta}}{C_{\delta}} + 
\tfrac{ (C_z^{H_1})^2}{c^{\bfB}C_{\delta}} \abs{\Omega} T
+ \bigl( \tfrac{9 }{c^{\bfA}} {+} \tfrac{2}{c^{\bfB}} \bigr) 
\tfrac{\abs{\Omega}^{\frac{1}{2}}}{C_{\delta}}T \\&+ 
\tfrac{1}{c^{\bfA} \cK C_{\delta} } T \norm[\L^{\infty} (0,T; \L^2(\Omega))]{f}^2,
\end{aligned}
\end{equation*}
and the conclusion follows with Gr\"onwall's lemma.
\eproof

Now we rewrite \eqref{eq:ent_eq_2_2} as follows 
\begin{equation*}
\dot z-\alpha \bfB^{-1}\Delta z=\bfB^{-1}f^z,
\end{equation*}
with
\(f^z\eqldef \tilde \bfQ^{\tra} \bfE (\ee(u){-} \bfQ \ze) -\partial_z H_1(z)-
\tilde \theta\partial_zH_2(z)-\psi\) and \(\psi \in\partial\Psi (\dot z)\). 
With assumption \eqref{eq:Psi.bdd} we infer that $\psi \in \L^{\infty}
(0, \tau; \L^{\infty}(\Omega))$ with 
$\norm[\L^{\infty}(\Omega)]{\psi(\cdot, t) } \le C^{\Psi}$ 
almost every $t \in (0, \tau)$. Furthermore, we can estimate $f^z$ as 
\begin{equation*}
\abs{f^z} \leq \norm{\tilde \bfQ} \norm{\bfE} \bigl(\abs{\ee(u)} {+} 
\norm{\tilde \bfQ} \abs{z} {+} \abs{\mathrm{Q}} \bigr) 
{+}( C_z^{H_1} {+} C_z^{H_2} 
\abs{\tilde \theta} )(1{+} \abs{z}) + C^{\Psi}.
\end{equation*}
Thus, using Lemma \ref{lemma3}, we infer first an estimate of $f^z$ in 
$\L^{\infty}(0,\tau; \L^2(\Omega))$
given by 
\begin{equation} 
\label{eqlp:6}
\begin{aligned}
&\norm[\L^{\infty}(0,\tau; \L^2(\Omega))]{f^z} 
\\&\le C 
\bigl(C\bigl(\norm[\W^{1,2}(\Omega)]{u^0},\norm[\W^{1,2}(\Omega)]{z^0}\bigl)
(X{+}1) 
\exp(c_0(X{+}1)\tau){+}X{+}1\bigr),
\end{aligned}
\end{equation}
where $C$ is a  constant independent of the initial data and $\tau$.
Hence, for any $q>8$, we have 
\begin{equation} \label{eqlp:7}
\norm[\L^q(0, \tau; \W^{2,2}(\Omega))]{z} \le C_z^q (\tau) 
\bigl( \norm[\L^{\infty}(0,\tau; \L^2(\Omega))]{f^z}  {+} 
\norm[\W^{2,2} (\Omega)]{z^0} \bigr),
\end{equation}
with a  non decreasing positive mapping $\tau \mapsto C_z^q( \tau)$ 
(see \cite{HieReh08QPMB,PruSch01SMRP}).
It follows that 
\begin{equation} \label{eqlp:8}
\begin{aligned}
& \norm[\L^q(0, \tau; \L^4(\Omega))]{f^z}\leq  \norm{\tilde \bfQ} 
\norm[\L^{\infty}(\Omega)]{\bfE} \bigl(\norm[\L^q(0, \tau; 
\L^4(\Omega))]{\ee(u)} \\& 
{+} C_1 \tau^{\frac{1}{q}} \norm{\tilde \bfQ} \norm[\L^{\infty}(0, \tau; 
\W^{1,2} (\Omega))]{z}{+}  \tau^{\frac{1}{q}} \abs{\Omega}^{\frac{1}{4}} 
\abs{\mathrm{Q}} \bigr)\\& 
 + C_z^{H_1}  \tau^{\frac{1}{q}} \bigl( C_1 \norm[\L^{\infty}(0, \tau; 
\W^{1,2} (\Omega))]{z}
{+}\abs{\Omega}^{\frac{1}{4}} \bigr) \\&+  C_z^{H_2} \| \tilde 
\theta\|_{\C^0([0, \tau]; \L^4(\Omega))} 
 \bigl( \tau^{\frac{1}{q}}{+}C_2 \norm[\L^{q}(0,\tau;\W^{2,2}(\Omega))]{z}
 \bigr) 
+ C^{\Psi} \tau^{\frac{1}{q}} \abs{\Omega}^{\frac{1}{4}},
\end{aligned}
\end{equation}
where $C_1$ and $C_2$ are the two generic constants involved in 
the continuous embeddings of $\W^{1,2}(\Omega)$ into $\L^4(\Omega)$ 
and $\W^{2,2} (\Omega)$ into $\L^{\infty} (\Omega)$, respectively. 
By combining Lemma \ref{lemma2} and Lemma \ref{lemma3}, we have 
\begin{equation*}
\begin{aligned}
& \|u\|_{\C^0([0, \tau]; \W^{1,4}(\Omega))} + \norm[\L^q(0,\tau;
\W^{1,4} (\Omega))]{\dot u} 
\\&\leq
 C^q_{u} (\tau) \Bigl(\tfrac{\tau^{\frac{2}{q}}}{2} C 
\bigl(\norm[\W^{1,2}(\Omega)]{u^0},  \norm[\W^{1,2}(\Omega)]{z^0}\bigr) 
(X{+}1) \exp(c_0 (X{+}1) \tau) \\&{+} \tfrac{\tau^{\frac{2}{q}}}{2} X  {+} 
\norm[\W^{1,4}(\Omega)]{u^0} {+}2 \Bigr),
\end{aligned}
\end{equation*}
and gathering \eqref{eqlp:6}, \eqref{eqlp:7} and  \eqref{eqlp:8}, 
we infer that
\begin{equation*}
\norm[\L^q(0, \tau; \L^4(\Omega))]{f^z} \le 
C_{f^z}^q \bigl(\tau, \norm[\W^{1,4} (\Omega)]{u^0}, 
\norm[\W^{2,2}(\Omega)]{z^0}\bigr)  
(X{+}1)^2 \exp( c_0 (X{+}1) \tau), 
\end{equation*}
where $  C_{f^z}^q \bigl(\tau, \norm[\W^{1,4} (\Omega)]{u^0}, 
\norm[\W^{2,2}(\Omega)]{z^0}\bigr) $ is a non decreasing positive  
function  of each of its arguments.

Using classical maximal regularity results for parabolic equations 
(\cite{HieReh08QPMB,PruSch01SMRP}), we may obtain an analogous 
estimate for $\norm[\L^q(0, \tau; \L^4(\Omega))]{\dot z}$. 
More precisely, there exists $  C^q_{u, z} \bigl(\tau,
\norm[\W^{1,4}(\Omega)]{u^0}, 
\norm[\W^{2,4}(\Omega)]{z^0}\bigr) $, which  is a  
non decreasing positive function  of each of its arguments, such that 
\begin{equation*}
\norm[\L^q(0, \tau; \L^4(\Omega))]{\dot z} \le C^q_{u, z}  \bigl(\tau, 
\norm[\W^{1,4} (\Omega)]{u^0}, \norm[\W^{2,4} (\Omega)]{z^0}\bigr)  
(X{+}1)^2 \exp( c_0 (X{+}1) \tau), 
\end{equation*}
and 
\begin{equation*}
\norm[\L^q(0, \tau; \L^4(\Omega))]{\ee (\dot u)}  
\le C^q_{ u, z}  \bigl(\tau, \norm[\W^{1,4} (\Omega)]{u^0}, 
\norm[\W^{2,4} (\Omega)]{z^0}\bigr)(X{+}1) \exp( c_0 (X{+}1) \tau) . 
\end{equation*}
 Finally, we have
\begin{equation*}
\begin{aligned}
& \norm[\L^{q/2}(0, \tau; \L^2(\Omega))]{f^{\tilde \theta}}  
\le \norm{\bfA} \norm[\L^q(0, \tau; \L^4(\Omega))]{\ee (\dot u)}^2 
+ \norm{\bfB} \norm[\L^q(0, \tau; \L^4(\Omega))]{\dot z}^2 \\
& +    \| \tilde \theta\|_{\C^0([0, \tau]; \L^4(\Omega))} 
\bigl( 3 \beta \tau^{\frac{1}{q}}  \norm[\L^q(0, \tau; 
\L^4(\Omega))]{\ee (\dot u)}\\&
{+} C_z^{H_2} \bigl(\tau^{\frac{1}{q}} {+} C_2 \norm[\L^q(0, \tau; 
\W^{2,2} (\Omega))]{z} \bigr) 
\norm[\L^q(0, \tau; \L^4(\Omega))]{\dot z} \bigr) 
+ C^{\Psi} \abs{\Omega}^{\frac{1}{4}} \tau^{\frac{1}{q}} 
\norm[\L^q(0, \tau; \L^4(\Omega))]{\dot z} \\
& \le C_{f^{\tilde \theta}}^q  (\tau, \norm[\W^{1,4}(\Omega)]{u^0}, 
\norm[\W^{2,4} (\Omega)]{z^0}) 
(X {+}1)^4 \exp( 4 c_0 (X{+}1) \tau), 
\end{aligned}
\end{equation*}
where once again $ C_{f^{\tilde\theta}}^q \bigl(\tau, \norm[\W^{1,4}(\Omega)]{u^0}, 
\norm[\W^{2,4} (\Omega)]{z^0}\bigr) $ is a  non decreasing positive function  
of each of its arguments.
It follows that 
\begin{equation} 
\label{eqlp:final}
\begin{aligned}
& \| \theta\|_{\C^0([0, \tau]; \L^4(\Omega))}  =  
\| \Phi^{\tilde \theta, \theta}_{\tau} (\tilde \theta) 
\|_{\C^0([0, \tau]; \L^4(\Omega))}\\& 
\leq C_1 C_{\theta} \exp 
\bigl(\tfrac{\tau}{2}\bigr) \bigl( \norm[\W^{1,2}(\Omega)]{\theta^0} {+} 
\norm[\L^{2}(0, \tau; \L^2(\Omega))]{f^{\tilde \theta}} \bigr) \\&
\leq  C_1 C_{\theta} \exp\bigl(\tfrac{\tau}{c^c}\bigr) \bigl( \norm[\W^{1,2} 
(\Omega)]{\theta^0}
 {+} C_{f^{\tilde \theta}}^q  \bigl(\tau, \norm[\W^{1,4}(\Omega)]{u^0}, 
\norm[\W^{2,4}(\Omega)]{z^0}\bigr)\\& 
\times\tau^{\frac{q-4}{2q}} 
(X {+}1)^4 \exp( 4 c_0 (X{+}1) \tau) 
 \bigr) \\
& \le  C^q \bigl(\tau, \norm[\W^{1,4}(\Omega)]{u^0}, 
\norm[\W^{2,4}(\Omega)]{z^0}, 
\norm[\W^{1,2}(\Omega)]{\theta^0}\bigr) \\&  
\times(X{+}1)^4 \exp( 4 c_0 (X{+}1) 
\tau), 
\end{aligned}
\end{equation}
where 
\begin{equation*}
\begin{aligned}
&  C^q \bigl(\tau, \norm[\W^{1,4}(\Omega)]{u^0}, 
\norm[\W^{2,4}(\Omega)]{z^0}, \norm[\W^{1,2}(\Omega)]{\theta^0}\bigr) \\
& 
\eqldef C_1 C_{\theta} \exp \bigl( \tfrac{\tau}{c^c}\bigr) \bigl(  \norm[\W^{1,2} 
(\Omega)]{\theta^0}  
+  \tau^{\frac{q-4}{2q}}  C_{f^{\tilde \theta}}^q \bigl(\tau, \norm[\W^{1,4}
(\Omega)]{u^0}, \norm[\W^{2,4} (\Omega)]{z^0}\bigr) \bigr).
\end{aligned}
\end{equation*}
Let us fix now $q>8$ and define the mapping $\gamma^q$ by 
\begin{equation*}
\gamma^q : R^{\theta} \mapsto g^q((\beta^2 {+}(C_z^{H_2})^2 )  
(R^{\theta})^2) - R^{\theta},
\end{equation*}
with 
\begin{equation*}
g^q(X) \eqldef
 C^q  \bigl(T, \norm[\W^{1,4}(\Omega)]{u^0}, \norm[\W^{2,4}(\Omega)]{z^0}, 
\norm[\W^{1,2}(\Omega)]{\theta^0}\bigr) 
 (X {+}1)^4 \exp( 4 c_0 (X{+}1) T),
\end{equation*}
for all $X \ge 0$. 
Observing that $X \mapsto g^q(X)$ is a continuous function, 
we can check that for any $R^{\theta} > g^q(0)$, there exists 
$\varepsilon_{q} >0$ such that $\gamma^q (R^{\theta}) < 0$ if 
\begin{equation*}
0 < \beta^2  + (C_z^{H_2})^2 < \tfrac{ \varepsilon_q}{(R^{\theta})^2}.
\end{equation*} 
Assuming that this condition holds, \eqref{eqlp:final} shows 
that  ${\mathcal C} \eqldef {\bar B}_{\C^0([0, T]; 
\L^4(\Omega))} (0, R^{\theta})$ is a closed convex bounded  
subset of $\C^0([0, T]; \L^4(\Omega))$ such that 
$\Phi^{\tilde \theta, \theta}_{T} ({\mathcal C}) \subset {\mathcal C}$. 
By using once again Schauder's fixed point theorem we may conclude 
that problem \eqref{eq:ent_eq}--\eqref{eq:init_cond} admits a 
global solution $(u, z, \theta)$ on $[0,T]$.

%%%%%%%%%%%%%%%%%%%%%%%%%%%%%%%%%%%%%%%%%%%%%%%%%%%%%%%%%%%%%%%%%%%%%%%%%%%%%%%%

\section{Examples}
\label{sec:examples}

In this concluding section, we present two classes of materials 
which fit our modelization, namely visco-elasto-plastic materials and 
SMA undergoing thermal expansion.

Indeed, in the both cases, an internal variable $z$ belonging to a finite 
dimensional real vector space is introduced to describe the inelastic 
strain due to plasticity or to phase transitions via the relation
\begin{equation*}
\ee^{\textrm{\bf inel}} = \bfQ z
\end{equation*}
where $z \mapsto \bfQ z $ is an affine mapping. 
The Helmholtz free energy is given by
\begin{equation*}
\begin{aligned}
W(\ee(u), z, \theta)\eqldef &
\tfrac{1}{2} \bfE (\ee(u){-}\bfQ z){:} 
(\ee(u){-}\bfQ z)+ 
\tfrac{\alpha}{2} \abs{\nabla z}^2 
\\&+ H(z, \theta) 
- c(\theta\ln(\theta){-}\theta) 
+ \beta \bfI{:}\ee(u),
\end{aligned}
\end{equation*}
where $H(z, \theta)$ is a hardening functional that may depend 
on the temperature, $\beta \bfI$, with $\beta \ge 0$,  is the 
isotropic thermal expansion tensor and $\alpha \ge 0$ is a coefficient 
that measures non local interaction effects for the internal variable. 
As usual $\bfE$ denotes the elasticity tensor, 
$\ee(u) \eqldef \tfrac12(\nabla u{+}\nabla
u^{\tra})$ is the infinitesimal strain tensor,  
and $c$ and $\kappa$ are the heat capacity and conductivity.

For visco-elasto-plastic models $\bfQ$ is linear, $H$ does not 
depend on $\theta$ and $\alpha =0$  while $\bfQ$ may be linear or 
affine as well,  $\alpha >0$ and $H$ depends on $\theta$ for 
SMA. Thus, by replacing $H (z, \theta) $ by an affine 
approximation $H_1(z) + \theta H_2(z)$, we may split $W(\ee(u), z, \theta)$ as
\begin{equation*} 
 W^{\textrm{mech}} (\ee(u), z) - W^{\theta} (\theta) + 
\theta  W^{\textrm{coup}} (\ee(u), z) 
\end{equation*}
with
\begin{equation*} 
\begin{aligned}
& 
W^{\textrm{mech}}(\ee(u), z)\eqldef 
\tfrac{1}{2} \bfE (\ee(u){-}\bfQ z){:} (\ee(u){-}\bfQ z) + H_1(z) 
+\tfrac{\alpha}2\abs{\nabla z}^2, \\
& W^{\theta} (\theta)\eqldef c (\theta\ln(\theta){-} \theta), \\
& W^{\textrm{coup}}(\ee(u), z) \eqldef \beta \bfI{:}\ee(u) + H_2(z) .
\end{aligned}
\end{equation*}
Let us illustrate this general setting with more precise modelizations. 
In the case of thermo-visco-elasto-plasticity, we can consider the 
Melan-Prager model corresponding to a linear kinematic hardening, i.e.
we have 
\begin{equation*}
H(z, \theta) \eqldef H_1(z) = \tfrac{1}{2} \bfL z {.} z 
\quad{\textrm{and}}\quad H_2(z) \equiv 0,
\end{equation*}
with a symmetric positive definite tensor 
$\bfL \in {\mathcal L}( {\mathcal Z}, {\mathcal Z})$, 
or the Prandtl-Reuss model for which 
$H(z, \theta) \equiv 0 = H_1(z) = H_2(z)$ (see \cite{Mau92TMPF}).

In the case of SMA, we can consider the  3D  macroscopic  
phenomenological model
introduced by Souza, Auricchio  et al.
(\cite{SoMaZo98TDMS, AurPet02IACR, AurPet04STPT}, or so-called mixture models
(see \cite{Miel00EMFM,HalGov02ARFE,GoMiHa02FEMV,MiThLe02VFRI,
GoHaHe07UBFE}). 
In the former case, 
$z \in {\mathcal Z} \eqldef\Erdev =\{ z \in \Ers: \ \bfI{:}z =0 \}$ and 
$\ee^{\textrm{\bf inel}} = \bfQ z = z$. Moreover the hardening functional 
is given by
\begin{equation*}
H_{\textrm{SA}}(z, \theta)\eqldef c_1 (\theta) \abs{z}+
c_2 (\theta) \abs{z}^2+\chi(z),
\end{equation*} 
where  \(\chi\) 
is the indicator function of the
ball \(\{z\in\Er^{3\times 3}_{\text{dev}}:\ \abs{z}\leq c_{3} (\theta)\}\). 
This coefficient \(c_3 (\theta) \) corresponds to  the maximum modulus 
of transformation
strain that can be obtained by alignment of martensitic variants while 
\(c_1 (\theta) >0\) is an activation threshold for
initiation of martensitic phase transformations and  \(c_2 (\theta)\) 
measures the
occurrence of hardening with respect to the internal variable
\(z\). 

In order to fit our regularity assumptions for the hardening functionals, 
which were assumed to be of class $\C^2$, we consider the regularization 
of $H_{\textrm{SA}}$ given by 
\begin{equation*} 
H_{\textrm{SA}}^{\delta}(z,\theta)\eqldef
c_1 (\theta) \sqrt{\delta^2{+}\abs{z}^2}
+ c_2 (\theta )\abs{z}^2+
\tfrac{((\abs{z}{-}c_3 (\theta))_+)^4}{\delta (1{+}\abs{z}^2)}, 
\end{equation*}
with \(0 < \delta\ll 1\),
(see also \cite{MiePet07TDPT} for another regularization of $H_{\textrm{SA}}$).

In the latter case, i.e. in so called mixture models, $z \in {\mathcal Z} 
\eqldef \Er^{N-1}  $ where $N \ge 2$ is the total number of phases and 
$\ee^{\textrm{\bf inel}} = \bfQ z $ is 
the effective transformation strain of the mixture, given by
\begin{equation*}
\bfQ z \eqldef \sum_{k=1}^{N-1} z_k \ee_k + \Bigl(1{-}\sum_{k=1}^{N-1} 
z_k\Bigr) \ee_N,
\end{equation*}
where $\ee_k$ is the transformation strain of the phase $k$. Then $z_1, 
\dots, z_{N-1}$
and $z_N \eqldef 1-\sum_{k=1}^{N-1} z_k$ can be interpreted as phase fractions 
and 
\begin{equation*}
H_{\textrm{mixt}}(z, \theta) = w(z, \theta) + \chi(z)
\end{equation*}
where $\chi$ is the indicator function of the set $[0,1]^{N-1}$. Once 
again we may consider a regularization of $H_{\textrm{mixt}}$ given by 
\begin{equation*}
H^{\delta}_{\textrm{mixt}} (z,\theta) = 
w(z, \theta) + \sum_{k=1}^{N-1} 
\tfrac{((-z_k)_+)^4{+} 
((z_k{-}1)_+)^4}{\delta (1{+}\abs{z_k}^2)},
\end{equation*}
with \(0< \delta\ll 1\).

%\begin{thebibliography}{SK}
\renewcommand{\arraystretch}{0.91}\small 

\bibliographystyle{my_alpha}
\bibliography{bib_PaoPet11}

%% Use the widest label as parameter.

%% Reference items may be numbered or have labels of your choice.
%% The author's surname PRECEDES the initial of the first name
%% The issue number is only given when the issues are paginated separately.
%% In book titles, first letters are capitalized.
%% Only journal volume numbers are boldfaced.

%%%%%%%%%%% To ease editing, use normal size:

%\normalsize
%\baselineskip=17pt

%%%%%%%%%%%%%

%\end{thebibliography}

\renewcommand{\arraystretch}{0.91}

\end{document}